
\magnification=1200

\hsize=14cm
\vsize=20.5cm
\parindent=0cm   \parskip=0pt
\pageno=1

\def\ind{\hskip 1cm\relax}


\hoffset=15mm    
\voffset=1cm    


\ifnum\mag=\magstep1
\hoffset=-0.5cm   
\voffset=-0.5cm   
\fi


\pretolerance=500 \tolerance=1000  \brokenpenalty=5000

\catcode`\@=11

\font\eightrm=cmr8         \font\eighti=cmmi8
\font\eightsy=cmsy8        \font\eightbf=cmbx8
\font\eighttt=cmtt8        \font\eightit=cmti8
\font\eightsl=cmsl8        \font\sixrm=cmr6
\font\sixi=cmmi6           \font\sixsy=cmsy6
\font\sixbf=cmbx6


\font\tengoth=eufm10       \font\tenbboard=msbm10
\font\eightgoth=eufm8      \font\eightbboard=msbm8
\font\sevengoth=eufm7      \font\sevenbboard=msbm7
\font\sixgoth=eufm6        \font\fivegoth=eufm5

\font\tencyr=wncyr10       
\font\eightcyr=wncyr8      
\font\sevencyr=wncyr7      
\font\sixcyr=wncyr6

\skewchar\eighti='177 \skewchar\sixi='177
\skewchar\eightsy='60 \skewchar\sixsy='60


\newfam\gothfam           \newfam\bboardfam
\newfam\cyrfam

\def\tenpoint{%
  \textfont0=\tenrm \scriptfont0=\sevenrm \scriptscriptfont0=\fiverm
  \def\rm{\fam\z@\tenrm}%
  \textfont1=\teni  \scriptfont1=\seveni  \scriptscriptfont1=\fivei
  \def\oldstyle{\fam\@ne\teni}\let\old=\oldstyle
  \textfont2=\tensy \scriptfont2=\sevensy \scriptscriptfont2=\fivesy
  \textfont\gothfam=\tengoth \scriptfont\gothfam=\sevengoth
  \scriptscriptfont\gothfam=\fivegoth
  \def\goth{\fam\gothfam\tengoth}%
  \textfont\bboardfam=\tenbboard \scriptfont\bboardfam=\sevenbboard
  \scriptscriptfont\bboardfam=\sevenbboard
  \def\bb{\fam\bboardfam\tenbboard}%
 \textfont\cyrfam=\tencyr \scriptfont\cyrfam=\sevencyr
  \scriptscriptfont\cyrfam=\sixcyr
  \def\cyr{\fam\cyrfam\tencyr}%
  \textfont\itfam=\tenit
  \def\it{\fam\itfam\tenit}%
  \textfont\slfam=\tensl
  \def\sl{\fam\slfam\tensl}%
  \textfont\bffam=\tenbf \scriptfont\bffam=\sevenbf
  \scriptscriptfont\bffam=\fivebf
  \def\bf{\fam\bffam\tenbf}%
  \textfont\ttfam=\tentt
  \def\tt{\fam\ttfam\tentt}%
  \abovedisplayskip=12pt plus 3pt minus 9pt
  \belowdisplayskip=\abovedisplayskip
  \abovedisplayshortskip=0pt plus 3pt
  \belowdisplayshortskip=4pt plus 3pt
  \smallskipamount=3pt plus 1pt minus 1pt
  \medskipamount=6pt plus 2pt minus 2pt
  \bigskipamount=12pt plus 4pt minus 4pt
  \normalbaselineskip=12pt
  \setbox\strutbox=\hbox{\vrule height8.5pt depth3.5pt width0pt}%
  \let\bigf@nt=\tenrm       \let\smallf@nt=\sevenrm
  \normalbaselines\rm}

\def\eightpoint{%
  \textfont0=\eightrm \scriptfont0=\sixrm \scriptscriptfont0=\fiverm
  \def\rm{\fam\z@\eightrm}%
  \textfont1=\eighti  \scriptfont1=\sixi  \scriptscriptfont1=\fivei
  \def\oldstyle{\fam\@ne\eighti}\let\old=\oldstyle
  \textfont2=\eightsy \scriptfont2=\sixsy \scriptscriptfont2=\fivesy
  \textfont\gothfam=\eightgoth \scriptfont\gothfam=\sixgoth
  \scriptscriptfont\gothfam=\fivegoth
  \def\goth{\fam\gothfam\eightgoth}%
  \textfont\cyrfam=\eightcyr \scriptfont\cyrfam=\sixcyr
  \scriptscriptfont\cyrfam=\sixcyr
  \def\cyr{\fam\cyrfam\eightcyr}%
  \textfont\bboardfam=\eightbboard \scriptfont\bboardfam=\sevenbboard
  \scriptscriptfont\bboardfam=\sevenbboard
  \def\bb{\fam\bboardfam}%
  \textfont\itfam=\eightit
  \def\it{\fam\itfam\eightit}%
  \textfont\slfam=\eightsl
  \def\sl{\fam\slfam\eightsl}%
  \textfont\bffam=\eightbf \scriptfont\bffam=\sixbf
  \scriptscriptfont\bffam=\fivebf
  \def\bf{\fam\bffam\eightbf}%
  \textfont\ttfam=\eighttt
  \def\tt{\fam\ttfam\eighttt}%
  \abovedisplayskip=9pt plus 3pt minus 9pt
  \belowdisplayskip=\abovedisplayskip
  \abovedisplayshortskip=0pt plus 3pt
  \belowdisplayshortskip=3pt plus 3pt
  \smallskipamount=2pt plus 1pt minus 1pt
  \medskipamount=4pt plus 2pt minus 1pt
  \bigskipamount=9pt plus 3pt minus 3pt
  \normalbaselineskip=9pt
  \setbox\strutbox=\hbox{\vrule height7pt depth2pt width0pt}%
  \let\bigf@nt=\eightrm     \let\smallf@nt=\sixrm
  \normalbaselines\rm}

\tenpoint

\def\pc#1{\bigf@nt#1\smallf@nt}         \def\pd#1 {{\pc#1} }


\catcode`\;=\active
\def;{\relax\ifhmode\ifdim\lastskip>\z@\unskip\fi
\kern\fontdimen2  -1.2 \fontdimen3 \string;}

\catcode`\:=\active
\def:{\relax\ifhmode\ifdim\lastskip>\z@\unskip\fi\penalty\@M\ \fi\string:}

\catcode`\!=\active
\def!{\relax\ifhmode\ifdim\lastskip>\z@
\unskip\fi\kern\fontdimen2  -1.1 \fontdimen3 \string!}

\catcode`\?=\active
\def?{\relax\ifhmode\ifdim\lastskip>\z@
\unskip\fi\kern\fontdimen2  -1.1 \fontdimen3 \string?}

\def\^#1{\if#1i{\accent"5E\i}\else{\accent"5E #1}\fi}
\def\"#1{\if#1i{\accent"7F\i}\else{\accent"7F #1}\fi}

\frenchspacing


\newtoks\auteurcourant      \auteurcourant={\hfil}
\newtoks\titrecourant       \titrecourant={\hfil}

\newtoks\hautpagetitre      \hautpagetitre={\hfil}
\newtoks\baspagetitre       \baspagetitre={\hfil}

\newtoks\hautpagegauche
\hautpagegauche={\eightpoint\rlap{\folio}\hfil\the\auteurcourant\hfil}
\newtoks\hautpagedroite
\hautpagedroite={\eightpoint\hfil\the\titrecourant\hfil\llap{\folio}}

\newtoks\baspagegauche      \baspagegauche={\hfil}
\newtoks\baspagedroite      \baspagedroite={\hfil}

\newif\ifpagetitre          \pagetitretrue





\def\raggedbottom{\topskip 10pt plus 36pt\r@ggedbottomtrue}



\def\pointir{\unskip . --- \ignorespaces}

\def\Bigbreak{\vskip-\lastskip\bigbreak}
\def\Medbreak{\vskip-\lastskip\medbreak}


\def\ctexte#1\endctexte{%
  \hbox{$\vcenter{\halign{\hfill##\hfill\crcr#1\crcr}}$}}


\long\def\ctitre#1\endctitre{%
    \ifdim\lastskip<24pt\vskip-\lastskip\bigbreak\bigbreak\fi
  		\vbox{\parindent=0pt\leftskip=0pt plus 1fill
          \rightskip=\leftskip
          \parfillskip=0pt\bf#1\par}
    \bigskip\nobreak}

\long\def\section#1\endsection{%
\vskip 0pt plus 3\normalbaselineskip
\penalty-250
\vskip 0pt plus -3\normalbaselineskip
\Bigbreak
\message{[section \string: #1]}{\bf#1\unskip}\pointir}

\long\def\sectiona#1\endsection{%
\vskip 0pt plus 3\normalbaselineskip
\penalty-250
\vskip 0pt plus -3\normalbaselineskip
\Bigbreak
\message{[sectiona \string: #1]}%
{\bf#1}\medskip\nobreak}

\long\def\subsection#1\endsubsection{%
\Medbreak
{\it#1\unskip}\pointir}

\long\def\subsectiona#1\endsubsection{%
\Medbreak
{\it#1}\par\nobreak}

\def\rem#1\endrem{%
\Medbreak
{\it#1\unskip} : }

\def\remp#1\endrem{%
\Medbreak
{\pc #1\unskip}\pointir}

\def\rema#1\endrem{%
\Medbreak
{\it #1}\par\nobreak}

\def\newparwithcolon#1\endnewparwithcolon{
\Medbreak
{#1\unskip} : }

\def\newparwithpointir#1\endnewparwithpointir{
\Medbreak
{#1\unskip}\pointir}

\def\newpara#1\endnewpar{
\Medbreak
{#1\unskip}\smallskip\nobreak}

\let\+=\tabalign

\def\signature#1\endsignature{\vskip 15mm minus 5mm\rightline{\vtop{#1}}}

\mathcode`A="7041 \mathcode`B="7042 \mathcode`C="7043 \mathcode`D="7044
\mathcode`E="7045 \mathcode`F="7046 \mathcode`G="7047 \mathcode`H="7048
\mathcode`I="7049 \mathcode`J="704A \mathcode`K="704B \mathcode`L="704C
\mathcode`M="704D \mathcode`N="704E \mathcode`O="704F \mathcode`P="7050
\mathcode`Q="7051 \mathcode`R="7052 \mathcode`S="7053 \mathcode`T="7054
\mathcode`U="7055 \mathcode`V="7056 \mathcode`W="7057 \mathcode`X="7058
\mathcode`Y="7059 \mathcode`Z="705A

\def\spacedmath#1{\def\packedmath##1${\bgroup\mathsurround=0pt ##1\egroup$}%
\mathsurround#1 \everymath={\packedmath}\everydisplay={\mathsurround=0pt }}

\def\nospacedmath{\mathsurround=0pt \everymath={}\everydisplay={} }


\long\def\th#1 #2\enonce#3\endth{%
   \Medbreak
   {\pc#1} {#2\unskip}\pointir{\it #3}\medskip}

\long\def\tha#1 #2\enonce#3\endth{%
   \Medbreak
   {\pc#1} {#2\unskip}\par\nobreak{\it #3}\medskip}


\long\def\Th#1 #2 #3\enonce#4\endth{%
   \Medbreak
   #1 {\pc#2} {#3\unskip}\pointir{\it #4}\medskip}

\long\def\Tha#1 #2 #3\enonce#4\endth{%
   \Medbreak
   #1 {\pc#2} #3\par\nobreak{\it #4}\medskip}


\def\decale#1{\smallbreak\hskip 28pt\llap{#1}\kern 5pt}
\def\decaledecale#1{\smallbreak\hskip 34pt\llap{#1}\kern 5pt}
\def\puce{\smallbreak\hskip 6pt{$\scriptstyle\bullet$}\kern 5pt}



\def\displaylinesno#1{\displ@y\halign{
\hbox to\displaywidth{$\@lign\hfil\displaystyle##\hfil$}&
\llap{$##$}\crcr#1\crcr}}


\def\ldisplaylinesno#1{\displ@y\halign{
\hbox to\displaywidth{$\@lign\hfil\displaystyle##\hfil$}&
\kern-\displaywidth\rlap{$##$}\tabskip\displaywidth\crcr#1\crcr}}


\def\eqalign#1{\null\,\vcenter{\openup\jot\m@th\ialign{
\strut\hfil$\displaystyle{##}$&$\displaystyle{{}##}$\hfil
&&\quad\strut\hfil$\displaystyle{##}$&$\displaystyle{{}##}$\hfil
\crcr#1\crcr}}\,}


\def\system#1{\left\{\null\,\vcenter{\openup1\jot\m@th
\ialign{\strut$##$&\hfil$##$&$##$\hfil&&
        \enskip$##$\enskip&\hfil$##$&$##$\hfil\crcr#1\crcr}}\right.}


\let\@ldmessage=\message

\def\message#1{{\def\pc{\string\pc\space}%
                \def\'{\string'}\def\`{\string`}%
                \def\^{\string^}\def\"{\string"}%
                \@ldmessage{#1}}}

\def\diagram#1{\def\normalbaselines{\baselineskip=0pt
\lineskip=10pt\lineskiplimit=1pt}   \matrix{#1}}



\def\up#1{\raise 1ex\hbox{\smallf@nt#1}}


\def\cf{{\it cf}} 

\def\qed{\raise -2pt\hbox{\vrule\vbox to 10pt{\hrule width 4pt
                 \vfill\hrule}\vrule}}

\def\cqfd{\unskip\penalty 500\quad\vrule height 4pt depth 0pt
width 4pt\medbreak}

\def\virg{\raise .4ex\hbox{,}}   


\def\build#1_#2^#3{\mathrel{
\mathop{\kern 0pt#1}\limits_{#2}^{#3}}}


\def\boxit#1#2{%
\setbox1=\hbox{\kern#1{#2}\kern#1}%
\dimen1=\ht1 \advance\dimen1 by #1 \dimen2=\dp1 \advance\dimen2 by #1
\setbox1=\hbox{\vrule height\dimen1 depth\dimen2\box1\vrule}%
\setbox1=\vbox{\hrule\box1\hrule}%
\advance\dimen1 by .4pt \ht1=\dimen1
\advance\dimen2 by .4pt \dp1=\dimen2  \box1\relax}

\def\date{\the\day\ \ifcase\month\or janvier\or f\'evrier\or mars\or
avril\or mai\or juin\or juillet\or ao\^ut\or septembre\or octobre\or
novembre\or d\'ecembre\fi \ {\old \the\year}}

\def\dateam{\ifcase\month\or January\or February\or March\or
April\or May\or June\or July\or August\or September\or October\or
November\or December\fi \ \the\day ,\ \the\year}

\def\moins{\mathop{\hbox{\vrule height 3pt depth -2pt
width 5pt}}}
\def\crog{{\vrule height 2.57mm depth 0.85mm width 0.3mm}\kern -0.36mm
[}

\def\crod{]\kern -0.4mm{\vrule height 2.57mm depth 0.85mm
width 0.3 mm}}
\def\moins{\mathop{\hbox{\vrule height 3pt depth -2pt
width 5pt}}}

\def\rond{\kern 1pt{\scriptstyle\circ}\kern 1pt}

\def\diagram#1{\def\normalbaselines{\baselineskip=0pt
\lineskip=10pt\lineskiplimit=1pt}   \matrix{#1}}

\def\hfl#1#2{\nospacedmath\smash{\mathop{\hbox to
12mm{\rightarrowfill}}\limits^{\scriptstyle#1}_{\scriptstyle#2}}}

\def\ghfl#1#2{\nospacedmath\smash{\mathop{\hbox to
25mm{\rightarrowfill}}\limits^{\scriptstyle#1}_{\scriptstyle#2}}}

\def\phfl#1#2{\nospacedmath\smash{\mathop{\hbox to
8mm{\rightarrowfill}}\limits^{\scriptstyle#1}_{\scriptstyle#2}}}

\def\pa{\S\kern.15em}

\def\Z{\hbox{\bf Z}}
\def\P{\hbox{\bf P}}

\def\C{\hbox{\bf C}}

\def\Aut{\mathop{\rm Aut}\nolimits}
\def\cad{c'est-\`a-dire}

\def\cf{{\it cf.\/}}

\def\codim{\mathop{\rm codim}\nolimits}

\def\GL{\mathop{\rm GL}\nolimits}

\def\Id{\hbox{\rm Id}}

\def\indp{\par\hskip0.5cm}
\def\isom{\simeq}

\def\loc{{\it loc.cit.\/}}
\def\long{\mathop{\rm long}\nolimits}
\def\lra{\longrightarrow}
\def\llra{\nospacedmath\hbox to 10mm{\rightarrowfill}}
\def\lllra{\nospacedmath\hbox to 15mm{\rightarrowfill}}

\def\no{n\up{o}}

\def\PGL{\mathop{\rm PGL}\nolimits}

\def\ra{\rightarrow}

\def\theo{th\'eor\`eme}

\def\tx{\kern -1.5pt -}

\def\A#1{{\cal A}_{#1}}

\def\cc#1{\hfill\kern .7em#1\kern .7em\hfill}

\def\dra{\ra\kern -3mm\ra}
\def\ldra{\lra\kern -3mm\ra}

\def\og{\leavevmode\raise.3ex\hbox{$\scriptscriptstyle\langle\!\langle$}}
\def\fg{\leavevmode\raise.3ex\hbox{$\scriptscriptstyle\,\rangle\!\rangle$}}

\catcode`\@=12

\showboxbreadth=-1  \showboxdepth=-1

\baselineskip=14pt
\spacedmath{2pt}
\font\eightrm=cmr10 at 8pt
\overfullrule=0mm
\parskip=1.7mm
\let\phi=\varphi
\let\epsilon=\varepsilon
\def\AutP{\Aut^0(\P)}
\def\tx{\hskip-.4mm-}
\input amssym.def
\input amssym
\def\vide{\varnothing}
\def\A{\hbox{\bf A}}

\null\smallskip
\ctitre
{\bf  TH\'EOR\`EMES DE CONNEXIT\'E POUR LES PRODUITS D'ESPACES
PROJECTIFS ET LES GRASSMANNIENNES } \endctitre
\smallskip
\centerline{Olivier {\pc DEBARRE}}
\vskip10mm
\ind Le \theo\ de connexit\'e de Fulton-Hansen \'enonce que si $X$ est une
vari\'et\'e irr\'eductible compl\`ete et $X\ra\P^n\times\P^n$ un morphisme
dont l'image est de codimension $<n$, l'image inverse de la diagonale est
{\it connexe} ([FH], [FL1]). Le point de d\'epart de cet article est une
remarque de Fulton et Lazarsfeld sugg\'erant que cette propri\'et\'e de la
diagonale pourrait \^etre essentiellement num\'erique. C'est ce que nous
v\'erifions dans la premi\`ere partie, en d\'emontrant une \theo\ de
connexit\'e analogue pour les morphismes \`a valeurs dans un produit
d'espaces projectifs (th. 2.2), qui entra\^ine en particulier que
l'\'enonc\'e ci-dessus reste valable si l'on remplace la diagonale de
$\P^n\times\P^n$ par n'importe quelle sous-vari\'et\'e de dimension $n$ qui
domine chaque facteur.

\ind Dans la seconde partie, nous \'etudions le m\^eme probl\`eme pour les
morphismes \`a valeurs dans une grassmannienne $G(d,\P^n)$. Dans [H],
Hansen montre que si $X$ est une vari\'et\'e irr\'eductible compl\`ete et
$f:X\ra G(d,\P^n)\times G(d,\P^n)$ un morphisme dont l'image est de
codimension $<n$, l'image inverse de la diagonale est connexe. Des exemples
montrent que cette borne d\'ecevante est la meilleure possible en
g\'en\'eral (\S 5). Notre but est  d'am\'eliorer ce r\'esultat en tenant
compte des propri\'et\'es num\'eriques de $f(X)$. Le \theo\ 7.1 montre que
{\it l'image inverse de la diagonale est connexe, pourvu qu'il existe des
partitions $\lambda=(\lambda_0,\ldots,\lambda_d)$ et
$\mu=(\mu_0,\ldots,\mu_d)$ v\'erifiant $\lambda_i+\mu_{d-i}<n-d$}
($i=0,\ldots,d$) {\it telles que}  $[f(X)]\cdot p_1^*\sigma_{\lambda}\cdot
p_2^*\sigma_{\mu}\ne 0$ ($\sigma_\lambda$ et $\sigma_\mu$ sont les classes
de Schubert). Ce r\'esultat contient celui de Hansen, mais n'est pas
optimal.

\ind Les r\'esultats de connexit\'e du type pr\'ec\'edent sont toujours
appliqu\'es dans la si\-tua\-tion suivante: on se donne des vari\'et\'es
irr\'eductibles compl\`etes $X$ et $Y$, des morphismes $f:X\ra  G(d,\P^n)$
et $g:Y\ra  G(d,\P^n)$, et l'on veut conclure \`a la connexit\'e de
$X\times_{G(d,{\bf P}^n)}Y$. Le \theo\ 7.1 mentionn\'e ci-dessus requiert
pour cela l'hypoth\`ese $[f(X)]\cdot [g(Y)]\cdot
(\sigma_{1,\ldots,1}+\sigma_{n-d})\ne 0$.  Dans \S 8, on montre qu'on peut
affaiblir cette hypoth\`ese lorsque $Y$ est une sous-vari\'et\'e de
Schubert de $G(d,\P^n)$, ou une intersection de vari\'et\'es de Schubert
sp\'eciales. Sans entrer dans les d\'etails, mentionnons simplement que
lorsque $Y$ est une vari\'et\'e de Schubert associ\'ee \`a une partition
$\mu$ telle que $\mu_0>\cdots >\mu_r>\mu_{r+1}=0$ (par exemple une
vari\'et\'e de Schubert sp\'eciale), la condition requise est que
l'intersection de  $f(X)$ avec chacune des classes de Schubert de type
$(\mu_i+1,\ldots,\mu_i+1,\mu_{i+1},\ldots,\mu_r)$ soit non nulle
($i=0,\ldots,r$).

\ind Un dernier mot enfin sur les m\'ethodes. Il existe essentiellement
trois  approches aux \theo s de connexit\'e. La plus ancienne,
 celle de Grothendieck ([G]), utilise la g\'eom\'etrie formelle; elle fut
ensuite reprise et d\'evelopp\'ee par Hironaka et Matsumura, Hartshorne,
Ogus, Speiser et Faltings ([F]). La seconde (celle de [FH]) consiste \`a
prouver d'abord un r\'esultat d'irr\'eductibilit\'e du type Bertini,
valable en toute caract\'eristique et sans aucune hypoth\`ese de
propret\'e. Lorsque  les vari\'et\'es sont compl\`etes, on passe \`a la
connexit\'e en utilisant la factorisation de Stein. La derni\`ere m\'ethode
fut vraisemblablement initi\'ee par Mumford dans une lettre \`a Fulton de
1978, et utilis\'ee par la suite par Sommese et Van de Ven, Nori ([N]) et
moi-m\^eme ([D1]); elle exploite l'existence d'une action transitive d'un
groupe alg\'ebrique et fournit directement des r\'esultats de connexit\'e;
il faut supposer les vari\'et\'es compl\`etes, et la caract\'eristique du
corps de base nulle.

\ind C'est la seconde m\'ethode qui est adopt\'ee ici; nos r\'esultats sont
donc d\'emontr\'es sur un corps alg\'ebriquement clos de caract\'eristique
quelconque, et tous les r\'esultats de connexit\'e mentionn\'es ci-dessus
ont des analogues de type Bertini (dont ils sont cons\'equences), valables
sans hypoth\`ese de propret\'e. En particulier, on montre dans \S 6 que si
$f:X\ra G(d,\P^n)$ est un morphisme et $H$ un hyperplan g\'en\'eral de
$\P^n$, $f^{-1}\bigl( G(d,H)\bigr)$ est irr\'eductible si $f(X)$ rencontre
$G(d,M)$ pour tout sous-espace g\'en\'eral $M$ de $\P^n$ de codimension $2$
(le \theo\ de Bertini usuel correspond au cas $d=0$). On passe de cet
\'enonc\'e \`a celui sur l'image inverse de la diagonale de
$G(d,\P^n)\times G(d,\P^n)$ par une astuce de Deligne. \bigskip {\bf
Notations et conventions}

\ind On adopte des conventions analogues \`a celles de [FL1]: les
\'enonc\'es se rapportent \`a un corps de base alg\'ebriquement clos
arbitraire $k$ (cadre \og alg\'ebrique\fg ), sauf ceux marqu\'es $(k=\C)$,
pour lesquels le corps de base est $\C$, et la topologie la topologie
usuelle (cadre \og topologique\fg ). Dans le cadre topologique, si $f:Y\ra
X$ est une application continue, avec $X$ connexe, on \'ecrira
$\pi_1(Y)\dra\pi_1(X)$ pour signifier qu'il existe $y\in Y$ tel que
l'homomorphisme induit $\pi_1(f):\pi_1(Y,y)\ra \pi_1\bigl( X,f(y)\bigr)$
soit surjectif. Lorsque $Y$ est connexe, cette propri\'et\'e ne d\'epend
pas du choix du point $y$.

\ind Si $f:X\ra S$ est un $S$\tx sch\'ema et $\gamma$ un automorphisme de
$S$, on notera $^\gamma X$ le\break $S$\tx sch\'ema $\gamma f:X\ra S$.

\ind Les sous-vari\'et\'es seront toujours {\it ferm\'ees} dans
l'espace ambiant. Connexe (resp. irr\'eductible) signifie connexe
(resp. irr\'eductible) et non vide.

\ind Soit $L$ un espace projectif; on notera $G(d,L)$ la grassmannienne des
sous-espaces lin\'eaires de $L$ de dimension $d$ et, pour tout $u\in
G(d,L)$, $\Lambda_u$ l'espace lin\'eaire correspondant \`a $u$.\vskip1cm

\bigskip \centerline {I. PRODUITS D'ESPACES PROJECTIFS} \bigskip \ind Dans
cette partie, on fixe des entiers positifs $n_1,\ldots,n_r$; on note $\P$
le produit $\P^{n_1}\times \cdots\times \P^{n_r}$ et $\AutP$ le groupe
$\prod_{i=1}^r\PGL(n_i+1,k)$ agissant diagonalement sur $\P$.

\ind Pour toute partie non vide $I$ de $\{ 1,\ldots,r\}$, on note
$n_I=\sum_{i\in I}n_i$, $\P_I=\prod_{i\in I}\P^{n_i}$ et $p_I$ la
projection $\P\ra\P_I$.  \bigskip

{\bf 1. Un \theo\ de Bertini}

\medskip

{\pc LEMME} 1.1.-- {\it Soient $Y$ une vari\'et\'e irr\'eductible
compl\`ete, $X$ une sous-vari\'ete irr\'eductible de $\P^n\times Y$, et
$p:X\ra\P^n$, $q:X\ra Y$ les deux projections. Soit $L$ un sous-espace
lin\'eaire de $\P^n$, g\'en\'eral de codimension $\le \dim p(X)$. Alors}
$$\dim q\bigl( p^{-1}(L)\bigr)=\min\bigl( \dim q(X),\dim (X)-\codim
(L)\bigr)\ .$$ {\bf D\'emonstration}. Compte tenu du \theo\ de Bertini
([FL1], th. 1.1), il suffit de traiter le cas o\`u $L$ est un hyperplan. On
a alors $\dim p(X)\ge 1$, de sorte que $p^{-1}(L)$ est un diviseur de $X$
qui varie sans point base. Soit $a$ la dimension d'une fibre g\'en\'erale
de la projection $X\ra q(X)$, c'est-\`a-dire $a= \dim(X)-\dim q(X)$.
Puisque $X$ est irr\'eductible, les diviseurs $D$ de $X$ tels que la
codimension de $q(D)$ dans $q(X)$ soit $>1$ sont en nombre fini. Si $a=0$,
on a donc, pour $L$ g\'en\'eral, $$\dim q\bigl( p^{-1}(L)\bigr)=\dim q(X)
-1=\dim (X) -1\ .$$ \ind Si $a>0$, l'hyperplan $L$ rencontre chaque fibre
de la projection propre $X\ra q(X)$, de sorte que  $q\bigl(
p^{-1}(L)\bigr)=q(X)$. Ceci prouve le lemme.\cqfd

\ind On aura aussi besoin du r\'esultat suivant, pour lequel je n'ai pas
trouv\'e de r\'ef\'erence ad\'equate.

{\pc LEMME} 1.2.-- {\it Soient $X$ une vari\'et\'e unibranche sur un corps
alg\'ebriquement clos de ca\-rac\-t\'eristique nulle, $f:X\ra\P^n$ un
morphisme, et $L$ un sous-espace lin\'eaire g\'en\'eral de $\P^n$. Alors
$f^{-1}(L)$ est unibranche.}

 {\bf D\'emonstration}. La normalisation $h:\tilde X\ra X$ \'etant un
hom\'eomorphisme, il en est de m\^eme de $(fh)^{-1}(L)\ra f^{-1}(L)$. Il
suffit donc de montrer que $(fh)^{-1}(L)$ est {\it normal}, et on se
ram\`ene ainsi au cas o\`u $X$ est normal. La d\'emonstration est alors
classique (\cf\ [J], [K], remarque (7)) et proc\`ede comme suit: posons
$Z=\{\,(x,u)\in X\times G(d,\P^n)\mid f(x)\in \Lambda_u\,\}$; la premi\`ere
projection $Z\ra X$ est lisse, de sorte que $Z$ est normal; par
cons\'equent, la fibre g\'en\'erique de la seconde projection $q:Z\ra
G(d,\P^n)$ est normale (ses anneaux locaux sont des anneaux locaux de $Z$),
donc g\'eom\'etriquement normale sur le corps $K\bigr( G(d,\P^n)\bigr)$
puisque ce dernier est de caract\'eristique nulle ([Gr1], prop. 6.7.7).
L'ensemble des points $u$ de $G(d,\P^n)$ tels que $q^{-1}(u)$ soit normal
\'etant localement constructible ([Gr2], prop. 9.9.4), il contient un
ouvert dense .\cqfd

\ind Montrons maintenant un \theo\ de Bertini pour les produits d'espaces
projectifs. \bigskip {\pc TH\'EOR\`EME} 1.3.-- {\it Supposons donn\'es une
vari\'et\'e irr\'eductible $X$, un morphisme $f:X\ra \P$, et, pour chaque
$i=1,\ldots,r$, un sous-espace lin\'eaire g\'en\'eral $L_i$ de
 $\P^{n_i}$, tels que\break  $\dim p_I\bigl(\overline{f(X)}\bigr)\ge
\sum_{i\in I}\codim(L_i)$ pour toute partie non vide $I$ de $\{
1,\ldots,r\}$; notons $L= L_1\times\cdots\times L_r$.

\indp {\rm 1)} $f^{-1}(L)$ est non vide de codimension
$\sum_{i=1}^r\codim(L_i)$ dans $X$. \smallskip \indp {\rm 2)} Posons $J=\{
i\in\{ 1,\ldots,r\}\mid L_i\ne \P^{n_i}\}$ et supposons de plus que pour
toute partie $I$ de $\{ 1,\ldots,r\}$ rencontrant $J$, on ait $\dim
p_I\bigl(\overline{f(X)}\bigr)> \sum_{i\in I}\codim(L_i)$.

\ind {\rm a)}  $f^{-1}(L)$ est irr\'eductible;

\ind {\rm b)} $(k=\C)$ si $X$ est localement irr\'eductible, $\pi_1\bigl(
f^{-1}(L)\bigr)\dra \pi_1(X)$.}

{\bf D\'emonstration}. Quitte \`a projeter sur $\P_J$, on peut supposer $J=
\{ 1,\ldots,r\}$.
 Supposons donc $\dim p_I\bigl(\overline{f(X)}\bigr)\ge\sum_{i\in
I}\codim(L_i)$ (resp. $>$) pour toute partie non vide $I$ de $\{
1,\ldots,r\}$ et montrons 1) (resp. 2)) par  r\'ecurrence sur $r$. Lorsque
$r=1$, 1) est trivial et 2) d\'ecoule de [FL1], th. 1.1. Supposons $r>1$.
Si, pour chaque $i=1,\ldots,r$, on a $\dim
p_i\bigl(\overline{f(X)}\bigr)=\codim (L_i)$, il ressort de l'hypoth\`ese
$\dim \bigl(\overline{f(X)}\bigr)\ge \sum_{i=1}^r\codim(L_i)$ que
$\overline{f(X)}=\prod_{i=1}^rp_i\bigl(\overline{f(X)}\bigr)$, auquel cas
1) est trivial. Supposons donc $\dim p_1\bigl(\overline{f(X)}\bigr)> \codim
(L_1)$. Soit $I$ une partie non vide de $\{ 2,\ldots,r\}$. Par [FL1], th.
1.1, $X'= (p_1f)^{-1}(L_1)$ est ferm\'e irr\'eductible de codimension
$\codim (L_1)$ dans $X$, et $\overline{f(X')}=\overline{f(X)}\cap
p_1^{-1}(L_1)$. D'autre part, le lemme 1.1 appliqu\'e \`a la
sous-vari\'et\'e $p_{\{ 1\}\cup I}\bigl(\overline{f(X)}\bigr)$ de
$\P_1\times\P_I$ donne  $$\dim p_I\bigl(\overline{f(X')}\bigr)=\min
\bigl(\dim p_I\bigl(\overline{f(X)}\bigr),\dim p_{\{ 1\}\cup
I}\bigl(\overline{f(X)}\bigr)-\codim (L_1)\bigr) \ge \sum_{i\in I}\codim
(L_i)$$ (resp. $>$). On peut donc appliquer l'hypoth\`ese de r\'ecurrence
au morphisme\break $f'=p_{\{2,\ldots,k\} }f:X'\ra \P_{\{2,\ldots,k\} }$; la
propri\'et\'e 1) (resp. 2)a)) en r\'esulte. Sous l'hypoth\`ese resp\'ee, et
si $X$ est unibranche, il en est de m\^eme de $X'$ (lemme 1.2); le cas
$r=1$ entra\^ine $\pi_1(X')\dra\pi_1(X)$ tandis que la propri\'et\'e 2)b)
pour $X'$ entra\^ine  $$\pi_1\bigl(  f^{-1}(L)\bigr)=\pi_1\bigl(
f'^{-1}(L_2\times\cdots\times L_r)\bigr)\dra\pi_1(X')\ .$$\ind La
propri\'et\'e 2)b) pour $X$ en d\'ecoule.\cqfd

\ind Passons maintenant au cas des espaces lin\'eaires quelconques. \bigskip
{\pc TH\'EOR\`EME} 1.4.-- {\it Supposons donn\'es une vari\'et\'e
irr\'eductible $X$, un morphisme $f:X\ra \P$, et, pour chaque
$i=1,\ldots,r$, un sous-espace lin\'eaire $L_i$ de
 $\P^{n_i}$, tels que\break  $\dim p_I\bigl(\overline{f(X)}\bigr)\ge
\sum_{i\in I}\codim(L_i)$ pour toute partie non vide $I$ de $\{
1,\ldots,r\}$; notons $L= L_1\times\cdots\times L_r$.

\indp {\rm 1)} Si $f$ est propre au-dessus d'un ouvert $V$ de $\P$ et que
$L$ est contenu dans $V$, $f^{-1}(L)$ est non vide. \smallskip \indp {\rm
2)} Posons $J=\{ i\in\{ 1,\ldots,r\}\mid L_i\ne \P^{n_i}\}$ et supposons de
plus que pour toute partie $I$ de $\{ 1,\ldots,r\}$ rencontrant $J$, on ait
$\dim p_I\bigl(\overline{f(X)}\bigr)> \sum_{i\in I}\codim(L_i)$.

\ind {\rm a)}  Si $f$ est propre au-dessus d'un ouvert $V$ de $\P$ et que
$L$ est contenu dans $V$, alors $f^{-1}(L)$ est connexe;

\ind {\rm b)} $(k=\C)$ si $X$ est localement irr\'eductible, pour tout
voisinage ouvert $U$ de $L$ dans $\P$, on a $\pi_1\bigl(
f^{-1}(U)\bigr)\dra \pi_1(X)$.}

{\bf D\'emonstration}. Comme plus haut, on peut supposer $J= \{
1,\ldots,r\}$. Montrons 2)b); tout voisinage $U$ de $L$ contient des
produits $M_1\times\cdots\times M_r$ auxquels le point 2)b) du th. 1.3
s'applique, d'o\`u 2)b). Supposons maintenant $f$ propre au-dessus d'un
ouvert $V$ de $\P$ et montrons 1) et 2)a). Pour chaque $i=1,\ldots,r$,
notons $G_i$ la grassmannienne des sous-espaces lin\'eaires de $\P^{n_i}$ de
m\^eme codimension que $L_i$ et $G=G_1\times\cdots\times G_r$. Soit $W$
l'ouvert de $G$ qui consiste en les $(u_1,\ldots,u_r)$ tels que
$\Lambda_{u_1}\times\cdots\times \Lambda_{u_r}\subset V$; notons $$Z=\{\
(x,u_1,\ldots,u_r)\in X\times W\mid f(x)\in \Lambda_{u_1}\times\cdots\times
\Lambda_{u_r}\ \}\ .$$  \ind La premi\`ere projection r\'ealise $Z$ comme
un ouvert dans un fibr\'e en produits de grassmanniennes au-dessus de $X$,
de sorte que $Z$ est irr\'eductible. Comme $f$ est propre au-dessus de $V$,
la projection $q:Z\ra V$ est aussi propre. Si $\dim
p_I\bigl(\overline{f(X)}\bigr)\ge \sum_{i\in I}\codim(L_i)$ (resp. $>$) pour
toute partie non vide $I$ de $\{ 1,\ldots,r\}$, le th. 1.3 entra\^ine que
les fibres g\'en\'erales de $q$ sont non vides (resp. irr\'eductibles). Il
s'ensuit que $q$ est surjective (resp. que les fibres de $q$ sont connexes,
par un argument classique utilisant la factorisation de Stein de $q$ (\cf\
[FL1], th. 2.1). Ceci termine la d\'emonstration.\cqfd  \medskip

{\it Remarque} 1.5.-- On aura besoin de la version plus fine de 2)b)
suivante: {\it pour tout $x\in f^{-1}(L)$, l'homomorphisme $\pi_1\bigl(
f^{-1}(U),x\bigr)\ra \pi_1(X,x)$ est surjectif}. Cela se d\'emontre comme
dans [FL1], remark 2.2. \bigskip {\bf 2. Th\'eor\`emes de connexit\'e}

\ind On note $\Delta$ la diagonale de $\P\times \P $ et, pour tout
$\gamma\in\AutP$, $^\gamma\!\Delta$ l'image de $\Delta$ par l'automorphisme
$(x,y)\mapsto (x,\gamma y)$ de $\P\times\P$. \smallskip {\pc LEMME} 2.1.--
{\it Soient $X$ une vari\'et\'e irr\'eductible et $f:X\ra \P\times \P$ un
morphisme.

\indp {\rm 1)} Si $X$ est compl\`ete et que $\dim (p_I\times p_I)f(X)\ge
n_I$ pour toute partie non vide $I$ de $\{ 1,\ldots,r\}$, $f^{-1}(\Delta)$
est non vide. \smallskip \indp {\rm 2)} On suppose que $\dim (p_I\times
p_I)f(X)> n_I$ pour toute partie non vide $I$ de $\{ 1,\ldots,r\}$.

\ind {\rm a)} Pour $\gamma\in\AutP$ g\'en\'eral,
$f^{-1}(\,\!^\gamma\!\Delta)$ est irr\'eductible;

\ind {\rm b)} si $X$ est compl\`ete, $f^{-1}(\Delta)$ est connexe;

\ind {\rm c)} $(k=\C)$  si $X$ est localement irr\'eductible et compl\`ete,
$\pi_1\bigl( f^{-1}(\Delta)\bigr)\dra \pi_1(X)$.}

{\bf D\'emonstration}. Comme dans
 [FL1], p. 39, on utilise une astuce de Deligne. Pour chaque
$i=1,\ldots,r$, on consid\`ere des coordonn\'ees homog\`enes
$[x^{(i)},y^{(i)}]$ sur $\P^{2n_i+1}$, o\`u $x^{(i)}$ et $y^{(i)}$ sont des
$(n_i+1)$\tx uplets d'\'el\'ements de $k$; on note $V_i$ l'ouvert de
$\P^{2n_i+1}$ d\'efini par $x^{(i)}\ne 0$ et $y^{(i)}\ne 0$, et on pose
$V=V_1\times\cdots\times V_r$. D\'efinissons un morphisme $\phi:V\ra
\P\times\P$ par la relation
$$\phi([x^{(1)},y^{(1)}],\ldots,[x^{(r)},y^{(r)}])=(
[x^{(1)}],\ldots,[x^{(r)}],[y^{(1)}],\ldots,[y^{(r)}])\ .$$ \ind Soit
$\tilde \gamma=(\tilde \gamma_1,\ldots,\tilde
\gamma_r)\in\prod_{i=1}^r\GL(n_i+1,k)$, soit $\gamma$ son image dans
$\AutP$ et soit $^{\tilde \gamma_i}L_i\subset V_i$ l'espace lin\'eaire
d\'efini par les \'equations $x^{(i)}=\tilde \gamma_i(y^{(i)})$. Les
$^{\tilde\gamma_i}L_i$ forment un ouvert dense dans $G(n_i,\P^{2n_i+1})$ et
$\phi$ induit un isomorphisme entre $^{\tilde \gamma}L=\, ^{\tilde
\gamma_1}L_1\times\cdots\times\,\! ^{\tilde\gamma_r}L_r$ et
$^\gamma\!\Delta$.

\ind Notons $X'=X\times_{\bf P}V$ et $f':X'\ra V$ le morphisme canonique;
$f'^{-1}(\,\! ^{\tilde \gamma}L)$ est isomorphe \`a $f^{-1}(\,\!
^\gamma\!\Delta)$ et $f'$ est propre lorsque $X$ est compl\`ete. Le point
2)a) est alors cons\'equence du th. 1.3.2)a) et les points 1) et 2)b) du
th. 1.4. Le point 2)c) se d\'eduit du point 2)b) du th. 1.4 et de la
remarque 1.5, comme dans [FL1], th. 3.1. et cor. 3.3.\cqfd \medskip

{\pc TH\'EOR\`EME} 2.2.-- {\it Soient $X$ et $Y$ des vari\'et\'es
irr\'eductibles et $f:X\ra\P$ et $g:Y\ra\P$ des morphismes.

\indp {\rm 1)} Si $X$ et $Y$ sont compl\`etes et que $\dim p_If(X)+\dim
p_Ig(Y)\ge n_I$ pour toute partie non vide $I$ de $\{ 1,\ldots,r\}$,
$X\times_{\bf P} Y$ est non vide. \smallskip \indp {\rm 2)} On suppose que
$\dim p_If(X)+\dim p_Ig(Y)> n_I$ pour toute partie non vide $I$ de $\{
1,\ldots,r\}$.

\ind {\rm a)} Pour $\gamma\in\AutP$ g\'en\'eral, $X\times_{\bf P}\,\!
^\gamma Y$ est irr\'eductible;

\ind {\rm b)} si  $X$ et $Y$ sont compl\`etes, $X\times_{\bf P} Y$ est
connexe;

\ind {\rm c)} $(k=\C)$ si $X$ et $Y$ sont localement irr\'eductibles et
compl\`etes,\break $\pi_1(X\times_{\bf P} Y)\dra \pi_1(X\times Y)$.}

{\bf D\'emonstration}. Appliquer le lemme 2.1 au morphisme $(f,g):X\times
Y\ra\P\times\P$.\cqfd
 \bigskip {\pc COROLLAIRE} 2.3.-- {\it Soient $X$ une vari\'et\'e
irr\'eductible compl\`ete, $f:X\ra\P$ un morphisme et $Y$ une
sous-vari\'et\'e irr\'eductible de $\P$ tels que, pour toute partie non
vide $I$ de $\{ 1,\ldots,r\}$, on ait $\dim p_If(X)+\dim p_I(Y)>n_I$.

\ind {\rm a)} $f^{-1}(Y)$ est connexe;

\ind {\rm b)} $(k=\C)$ si $X$ est localement irr\'eductible,
$\pi_1\bigl(f^{-1}(Y)\bigr)\dra \pi_1(X)$.}

{\bf D\'emonstration}. Le point a) se d\'eduit du th. 2.2.2)b); le point b)
se d\'eduit du th. 2.2.2)c) comme dans [FL1], cor. 4.3 (consid\'erer la
normalis\'ee de $Y$).\cqfd

\ind Le lecteur remarquera que le lemme 2.1.2) est cons\'equence du cor.
2.3. \medskip \ind On d\'emontre aussi par les m\^emes m\'ethodes des
r\'esultats analogues \`a ceux des cor. 5.2 et 5.3 de [FL1], comme par
exemple:

{\pc COROLLAIRE} 2.4.-- {\it Soit $X$ une sous-vari\'et\'e irr\'eductible de
$\P$ telle que, pour toute partie non vide $I$ de $\{ 1,\ldots,r\}$, on ait
$\dim p_I(X)>{1\over 2}n_I$.

\ind {\rm a)} $\pi_1^{\rm alg}(X)=0$;

\ind {\rm b)} $(k=\C)$  $X$ est simplement connexe.} \bigskip

\ind Il est peut-\^etre plus parlant d'interpr\'eter les hypoth\`eses des
corollaires ci-dessus en termes des classes des sous-vari\'et\'es qui
interviennent. Il existe ([FMSS] ou [Fu], ex. 8.3.7) une d\'ecomposition de
K\"unneth pour les groupes de Chow  $$A^m(\P)=\bigoplus_{{\bf
m}=(m_1,\ldots,m_r)\atop m_1+\cdots+m_r=m} A^{\bf m}(\P )\ ,$$ o\`u $A^{\bf
m}(\P)=A^{m_1}(\P^{n_1})\otimes\cdots\otimes A^{m_r}(\P^{n_r})$. Soit
 $X$ une sous-vari\'et\'e irr\'eductible de $\P$ de codimension $m$; la
composante de sa classe $[X]$ dans $A^{\bf m}(\P)\isom\Z$ est l'entier
positif $$[X]_{\bf m}=X\cdot H_1^{n_1-m_1}\cdot\ldots\cdot H_r^{n_r-m_r}\
,$$ o\`u, pour $i=1,\ldots,r$, $H_i$ est l'image inverse dans $\P$ d'un
hyperplan de $\P^{n_i}$. On v\'erifie que, pour tout entier $a$ et toute
partie non vide $I$ de $\{ 1,\ldots,r\}$, on a
 $$\codim p_I(X)\le a\qquad{\Longleftrightarrow}\qquad \exists\ {\bf m}\quad
\sum_{i\in I}m_i\le a\quad{\rm et}\quad [X]_{\bf m}\ne 0\ .\leqno (2.5)$$
\ind Les hypoth\`eses du th. 2.2 se r\'e\'ecrivent alors: pour toute partie
non vide $I$ de $\{ 1,\ldots,r\}$, il existe ${\bf m}$ et ${\bf m'}$ avec
$\sum_{i\in I}(m_i+m'_i)\le n_I$ (resp. $<n_I$), tels que $[f(X)]_{\bf
m}\ne 0$ et $[g(Y)]_{\bf m'}\ne 0$. \bigskip

\ind Le cor. 2.3 entra\^ine que certaines sous-vari\'et\'es d'un produit
d'espaces projectifs ont des propri\'et\'es de connexit\'e analogues \`a
celles de la petite diagonale de $(\P^n)^r$, telles qu'elles sont
expos\'ees dans [FL1].

\ind On dira qu'une sous-vari\'et\'e  $Z$ d'une vari\'et\'e $P$ est {\it
encombrante} si elle rencontre toute sous-vari\'et\'e de $P$ de dimension
$\ge\codim (Z)$. Pour qu'une sous-vari\'et\'e irr\'eductible $Z$ de $\P$
soit encombrante, il faut et il suffit que $[Z]_{\bf m}$ soit non nul d\`es
que $A^{\bf m}(\P )$ l'est (\cad\ lorsque $m_i\le n_i$ pour tout $i$). On a
aussi:

\medskip {\pc PROPOSITION} 2.6.-- {\it Pour qu'une sous-vari\'et\'e
irr\'eductible $Z$ de $\P$ soit encombrante, il faut et il suffit que pour
toute partie non vide $I$ de $\{ 1,\ldots,r\}$, on ait} $$\dim p_I(Z)=\min
\bigl( \dim(Z), n_I\bigr)\ .$$

{\bf D\'emonstration de la proposition}. Supposons que $\dim p_I(Z)=\min
\bigl( \dim(Z), n_I\bigr)$
 pour toute partie non vide $I$ de $\{ 1,\ldots,r\}$; soient $Y$ une
sous-vari\'et\'e irr\'eductible de $\P$ de dimension $\ge\codim (Z)$ et $I$
une partie non vide de $\{ 1,\ldots,r\}$. Si $\dim p_I(Z)= n_I$, alors
$\dim p_I(Z)+\dim p_I(Y)\ge n_I$. Sinon, on a $\dim p_I(Z)= \dim (Z)$; comme
$$\dim p_I(Y)\ge \dim (Y)-(\dim(\P)-n_I)=n_I-\codim (Y)\ge n_I-\dim(Z)\ ,$$
on a encore $\dim p_I(Z)+\dim p_I(Y)\ge n_I$. Il r\'esulte du th. 2.2.1)
que $Z$ rencontre $Y$.

\ind Supposons inversement $Z$ encombrante.  Soit $I$ une partie non vide
de $\{ 1,\ldots,r\}$ telle que $p_I(Z)\ne\P_I$, et soit $L$ un sous-espace
lin\'eaire de $\P_I$ de codimension $\dim p_I(Z)+1$, disjoint de $p_I(Z)$.
Alors $Z$ ne rencontre pas $p_I^{-1}(L)$, de sorte que $$\dim (Z)<\codim
p_I^{-1}(L) = \dim p_I(Z)+1\ .$$ \ind Ceci entra\^ine $\dim (Z)=\dim
p_I(Z)$ et termine la d\'emonstration.\cqfd

\ind On dira qu'une sous-vari\'et\'e  $Z$ de $\P$ est {\it bonne} si $\dim
p_i(Z)=\dim(Z)$ pour tout $i=1,\ldots,r$.

{\pc PROPOSITION} 2.7.-- {\it Soient $X$ une vari\'et\'e irr\'eductible
compl\`ete, $f:X\ra\P$ un morphisme et $Z$ une sous-vari\'et\'e
irr\'eductible de $\P$ telle que $\dim f(X)>\codim (Z)$.

\ind {\rm a)} Si $Z$ est encombrante et que $\dim p_if(X)>0$ pour tout
$i=1,\ldots,r$, $f^{-1}(Z)$ est connexe;

\ind {\rm b)} si $Z$ est bonne, $f^{-1}(Z)$ est connexe.}

{\bf D\'emonstration}. Soit $I$ une partie non vide de $\{ 1,\ldots,r\}$.
Sous les hypoth\`eses de a), on a soit $\dim (Z)> n_I$, auquel cas
 $\dim p_I(Z)=n_I$ et $\dim p_If(X)+\dim p_I(Z)>n_I$; soit $\dim (Z)\le
n_I$, auquel cas  $\dim p_I(Z)=\dim (Z)$ et  $$\eqalign{\dim p_If(X)+\dim
p_I(Z) &\ge\dim f(X) -\sum_{j\notin I} n_j+\dim (Z)\cr &> \codim (Z)
-\sum_{j\notin I} n_j+\dim (Z) = n_I\ .\cr} $$ \ind Dans chacun de ces deux
cas, on peut appliquer le cor. 2.3. Sous l'hypoth\`ese b), on est toujours
dans le deuxi\`eme cas.\cqfd

{\it Remarques} 2.8.-- 1)  La petite diagonale de $(\P^n)^r$ est bonne; on
retrouve donc les r\'esultats de connexit\'e de [FL1].

\ind 2) Sous les hypoth\`eses du corollaire, \^etre encombrante n'est pas
suffisant en g\'en\'eral pour assurer la connexit\'e de $f^{-1}(Z)$, comme
le montre l'exemple suivant. Soit $\Gamma$ une courbe irr\'eductible dans
$\P^1\times\P^1$ telle que la premi\`ere projection $\Gamma\ra\P^1$ soit de
degr\'e $>1$. Soient $s:\P^1\times\P^r\ra\P^{2r+1}$ le plongement de Segre
et $\P=\P^1\times\P^{2r+1}$;
 l'image $Z$ de $\Gamma\times\P^r$ par $\Id\times
s:\P^1\times\P^1\times\P^r\ra \P$ est encombrante, mais, pour un point
g\'en\'eral $x$ de $\P^1$, $Z\cap p_1^{-1}(x)$ n'est pas connexe,  bien que
 $\dim (Z)+\dim p_1^{-1}(x)=r+1+2r+1>2r+2$. \bigskip {\bf 3. Classes des
sous-vari\'et\'es irr\'eductibles}

\ind On s'int\'eresse ici aux classes des sous-vari\'et\'es irr\'eductibles
$X$ de $\P$. On rappelle que l'on a not\'e $[X]_{\bf m}$ les composantes de
$[X]$ dans la d\'ecomposition de K\"unneth de $A^m(\P)$ (avec ${\bf
m}=(m_1,\ldots,m_r)$ et $m_1+\cdots+m_r=m=\codim(X)$).

\bigskip {\pc PROPOSITION} 3.1.-- {\it Soit $X$ une sous-vari\'et\'e
projective irr\'eductible de $ \P$. Pour tout $r$-uplet ${\bf
m}=(m_1,\ldots,m_r)$ d'entiers positifs de somme $\codim (X)$, la composante
$[X]_{\bf m}$ est non nulle si et seulement si, pour toute partie non vide
$I$ de $\{ 1,\ldots,r\}$, on a
 $$\sum_{i\in I}m_i\ge \codim p_I(X)\ .$$ \ind En particulier, l'ensemble
des ${\bf m}$ tels que $[X]_{\bf m}\ne 0$ est l'ensemble des points entiers
d'un ensemble convexe.}

{\bf D\'emonstration}. Soit ${\bf m}$ un  $r$-uplet; si les in\'egalit\'es
de la proposition sont v\'erifi\'ees, le th. 1.3.1) entra\^ine que $[X]_{\bf
m}$ est la classe d'un $0$-cycle effectif non nul. La r\'eciproque
r\'esulte de (2.5).\cqfd \medskip \ind L'in\'egalit\'e de Hodge force en
fait des conditions plus contraignantes sur les classes des
sous-vari\'et\'es irr\'eductibles de $\P$. Pour $\alpha\in\{ 1,\ldots,r\}$,
posons ${\bf
e}_{\alpha}=(\delta_{1,\alpha},\delta_{2,\alpha},\ldots,\delta_{r,\alpha})$;
si $\alpha,\beta\in\{ 1,\ldots,r\}$, on a alors
 $$[X]_{\bf m}^2\ge [X]_{{\bf m+e}_{\alpha}-{\bf e}_{\beta}}[X]_{{\bf m-
e}_{\alpha}+{\bf e}_{\beta}}\ .$$ \ind Ces in\'egalit\'es r\'esultent des
\theo s de Hodge et Bertini, et entra\^inent toute une s\'erie
d'in\'egalit\'es de convexit\'e connues sous le nom de \theo\ de
Teissier-Hovanski (\cf\ [T], [Ho], [G]), qui permettent en particulier de
retrouver la prop. 3.1.

\vskip1cm\centerline {II. GRASSMANNIENNES} \bigskip \ind On appelle {\it
partition} un $(d+1)$\tx uplet $\lambda=(\lambda_0,\ldots,\lambda_d)$
d'entiers tels que\break $n-d\ge\lambda_0\ge\cdots\ge\lambda_d\ge 0$; on
sous-entendra toujours $\lambda_i=0$ pour $i>d$, et on omettra m\^eme
souvent les $\lambda_i$ nuls. Notons $|\lambda
|=\lambda_0+\cdots+\lambda_d$ l'entier partitionn\'e; \`a $\lambda$ est
associ\'e un diagramme de Young, contenu dans un rectangle de $d+1$ lignes
(num\'erot\'ees de $0$ \`a $d$) et $n-d$ colonnes, et obtenu en pla\c cant
$\lambda_i$ bo\^ites dans la ligne $i$. Pour la partition $(4,3,2,2)$, cela
donne

$$\nospacedmath\vbox{\let\m=\multispan \offinterlineskip \halign{ &\hbox to
20pt{\vrule height 12pt depth 5pt width 0pt \hss$#$\hss}\cr
\noalign{ \hbox{\vrule height 0pt depth .4pt width 70.5pt \vrule height 0pt
depth 1pt width 52pt} \vskip -.5pt}
\m 1\vrule&\m 3
\hfill$\lambda_0$\hfill &\vrule width 1pt&&&\kern -17.5pt\vrule\cr
\noalign{ \hbox{\vrule height 0pt depth .4pt width 50.5pt
\vrule height 0pt depth 1pt width 21pt} \vskip -1pt}
\m 1\vrule &\m 2 \hfill $\lambda_1$\hfill &\vrule width 1pt&&&&\kern
-17.5pt\vrule\cr
 \noalign{ \hbox{\vrule height 0pt depth
.4pt width 30.5pt \vrule height 0pt depth 1pt width 21pt} \vskip -1pt}
 \m 1\vrule &\hfill \lambda_2 &\vrule width 1pt&&&&&\kern
-17.5pt\vrule\cr
 \noalign{ \hbox{\vrule height 0pt depth
.4pt width 30.5pt} \vskip -1pt}
 \m 1\vrule &\hfill
\lambda_3&\vrule width 1pt&&&&&\kern -17.5pt\vrule\cr
\noalign{ \hbox{\vrule height 1pt depth 0pt width 31.5pt} \vskip -1pt} \m
1\vrule width 1pt&&&&&&&\kern -17.5pt\vrule\cr
 \noalign{
\hbox{\vrule height .4pt width 122.4pt} } }}$$
 \ind On note $\bar\lambda$
la partition $(n-d-\lambda_d,\ldots,n-d-\lambda_0)$. Sur le diagramme de
Young de $\lambda$, elle se lit de bas en haut, \`a droite de l'escalier.

\ind On note aussi $\lambda^*$ la partition
$(\lambda^*_0,\ldots,\lambda^*_{n-d-1})$, d\'efinie par
$\lambda^*_i=\max\{\, j\mid \lambda_j>i\,\}$, pour $i=0,\ldots,n-d-1$. Sur
le diagramme de Young de $\lambda$, elle se lit verticalement, de gauche
\`a droite, au-dessus de l'escalier. Par exemple, $(4,3,2,2)^*=(4,4,2,1)$.
On notera que $\overline{\lambda^*}=\bar\lambda^*$.

\ind  Si $\lambda$ et $\mu$ sont deux partitions, on \'ecrit
$\lambda\le\mu$ si $\lambda_i\le\mu_i$ pour tout $i=0,\ldots,d$, et
$\lambda<\mu$ si $\lambda_i<\mu_i$ pour tout $i=0,\ldots,d$.

\ind  A toute partition $\lambda$ correspond une {\it classe de Schubert}
$\sigma_\lambda$ dans $A^{|\lambda|}\bigl( G(d,\P^n)\bigr)$. Pour tout
drapeau ${\bf L}=(L^{(0)}\subset \cdots\subset L^{(d)} \subset\P^n)$  avec
$\dim (L^{(i)})=n-d+i-\lambda_i$, $\sigma_\lambda$ est la classe de la
vari\'et\'e de Schubert $$\Sigma_{\bf L}=\{\ u\in G(d,\P^n)\mid \dim
\bigr(\Lambda_u\cap L^{(i)}\bigr)\ge i,\quad{\rm pour}\ \  i=0,\ldots,d\
\}\ .$$  \ind Pour $m\in\{ 0,\ldots,n-d\}$, la classe $\sigma_m$ est dite
{\it sp\'eciale}; les vari\'et\'es de Schubert associ\'ees sont
$$\Sigma_L=\{\ u\in G(d,\P^n)\mid \Lambda_u\cap L\ne\vide\ \}\ ,$$ o\`u
$\codim (L)=d+m$.

\ind Les classes de Schubert forment une base du $\Z$\tx module $A^*\bigl(
G(d,\P^n)\bigr)$. On a $\sigma_\lambda\cdot\sigma_{\bar\lambda}=1$ et
l'isomorphisme de dualit\'e $\phi:G(d,\P^n)\ra G\bigl(
n-d-1,(\P^n)^*\bigr)$ v\'erifie\break
$\phi^*(\sigma_{\lambda^*})=\sigma_\lambda$.

\ind Si $X$ est une sous-vari\'et\'e de $G(d,\P^n)$, sa classe dans
$A^*\bigl( G(d,\P^n)\bigr)$ s'\'ecrit $[X]=\sum_\lambda\
[X]_\lambda\,\sigma_\lambda$, avec
$[X]_\lambda=[X]\cdot\sigma_{\bar\lambda}$; les $[X]_\lambda$ sont des
entiers positifs non tous nuls. \bigskip
 {\bf 4. Calcul de Schubert}

\ind Voici deux r\'esultats \'el\'ementaires pour lesquels je n'ai pas
trouv\'e de r\'ef\'erence.
 \medskip  {\pc LEMME} 4.1.-- {\it On a $$\sigma_\lambda\cdot\sigma_\mu\ne
0\qquad{\Longleftrightarrow}\qquad \lambda\le\bar\mu\ .$$ \ind En
particulier, si $\sigma_\lambda\cdot\sigma_\mu\ne 0$ et
$\lambda'\le\lambda$, alors $\sigma_{\lambda'}\cdot\sigma_\mu\ne 0$.}

{\bf D\'emonstration}. Si la propri\'et\'e $\lambda\le\bar\mu$ n'est pas
v\'erifi\'ee, $\sigma_\lambda\cdot\sigma_\mu=0$ ([GH], p. 198).

\ind Supposons au contraire $\lambda\le\bar\mu$. Montrons par r\'ecurrence
sur $|\bar\mu|-|\lambda|$ que $\sigma_\lambda\cdot\sigma_\mu$ est non nul.
Si $|\lambda|=|\bar\mu|$, alors $\lambda=\bar\mu$ et
$\sigma_\lambda\cdot\sigma_\mu=1$. Si  $|\lambda|<|\bar\mu|$, on choisit
$i_0$ minimal tel que $\lambda_{i_0}+\mu_{d-i_0}<n-d$. La partition
$\lambda'$ d\'efinie par $\lambda'_i=\lambda_i+\delta_{i,i_0}$ v\'erifie
$\lambda'\le\bar\mu$ et $\lambda'_0\le n-d$, et l'hypoth\`ese de
r\'ecurrence entra\^ine $\sigma_{\lambda'}\cdot\sigma_\mu\ne 0$. Comme
$\lambda'$ intervient avec un coefficient non nul dans la d\'ecomposition
du produit $\sigma_1\cdot\sigma_\lambda$ en somme de classes de Schubert
donn\'ee par la formule de Pieri ([Fu], p. 271), on a aussi
$\sigma_1\cdot\sigma_\lambda\cdot\sigma_\mu\ne 0$, d'o\`u le lemme.\cqfd
\smallskip  \ind Le lemme 4.1 entra\^ine: {\it soient $X$ et $Y$ des
sous-vari\'et\'es de  $G(d,\P^n)$; pour que $X\cap Y$ soit non vide, il
faut et il suffit qu'il existe des partitions $\lambda$ et $\mu$ avec
$[X]_\lambda\ne 0$, $[Y]_\mu\ne 0$ et $\lambda\le\bar\mu$.}

\ind On aura aussi besoin du r\'esultat suivant sur les produits de classes
de Schubert sp\'eciales.
 \medskip  {\pc LEMME} 4.2.-- {\it Soient $\ell_0,\ldots,\ell_r$ des
entiers avec $n-d\ge\ell_0\ge\cdots\ge\ell_r\ge 0$, et $\lambda$ une
partition. Pour que
$\sigma_{\bar\lambda}\cdot\sigma_{\ell_0}\cdots\sigma_{\ell_r}$ soit non
nul, il faut et il suffit que
$$\ell_0+\cdots+\ell_i\le\lambda_0+\cdots+\lambda_i$$  pour tout $i=
0,\ldots,r$.}

\ind En prenant $|\lambda|=\sum_{i=0}^r\ell_i$, on obtient une description
explicite des classes de Schubert qui apparaissent avec un coefficient non
nul dans $\sigma_{\ell_0}\cdots\sigma_{\ell_r}$. Le lemme entra\^ine aussi
que {\it ce produit est non nul lorsque}
 $\sum_{i=0}^r\ell_i\le (d+1)(n-d)$.

{\bf D\'emonstration du lemme}. On proc\`ede par r\'ecurrence sur $r$.
Lorsque $r=-1$, il n'y a rien \`a d\'emontrer. Supposons l'\'equivalence
d\'emontr\'ee pour $r-1\ge -1$ et prenons des entiers
$\ell_0,\ldots,\ell_r$ tels que $n-d\ge\ell_0\ge\cdots\ge\ell_r\ge 0$.

\ind Supposons
$\sigma_{\bar\lambda}\cdot\sigma_{\ell_0}\cdots\sigma_{\ell_r}$ non nul; il
suffit par hypoth\`ese de r\'ecurrence de montrer l'in\'egalit\'e
$\ell_0+\cdots+\ell_r\le\lambda_0+\cdots+ \lambda_r$. Posons $s=\min (r,d)$
et soit $\bar\lambda'$ la partition $\bar\lambda$ tronqu\'ee apr\`es
$\bar\lambda_s$. On a alors $\sigma_{\bar\lambda'}
\cdot\sigma_{\ell_0}\cdots\sigma_{\ell_r}\ne 0$ (lemme 4.1). Ces cycles
vivent en fait dans $G(s,\P^{s+n-d})$, de sorte que
$$\bar\lambda_0+\cdots+\bar\lambda_s+\ell_0+\cdots+\ell_r\le \dim
G(s,\P^{s+n-d})=(s+1)(n-d)\ ,$$ ce qui d\'emontre l'in\'egalit\'e
cherch\'ee.

\ind Supposons inversement que les in\'egalit\'es du lemme soient
v\'erif\'ees. Soit $s\in\{ 0,\ldots,d\}$ l'entier tel que
$\lambda_{s+1}+\cdots+\lambda_d<\ell_r\le \lambda_s+\cdots+\lambda_d$ et
posons $\lambda'_i=\lambda_i$ pour $0\le i<s$,
$\lambda'_s=\lambda_s+\cdots+\lambda_d-\ell_r$ et $\lambda'_i=0$ pour
$i>s$. La formule de Pieri s'\'ecrit
$$\sigma_{\bar\lambda}\cdot\sigma_{\ell_r}=\sum_{|\mu|=|\lambda|-\ell_r\atop
\lambda_0\ge\mu_0\ge\cdots\ge\lambda_d\ge\mu_d\ge 0}\sigma_{\bar\mu}\ ,$$
de sorte que $\sigma_{\bar\lambda'}$ intervient avec un coefficient non nul
dans $\sigma_{\bar\lambda}\cdot\sigma_{\ell_r}$. Pour $ s\le i<r $, on a
$$\eqalign{\lambda'_0+\cdots+\lambda'_i&=\lambda_0+\cdots+
\lambda_{s-1}+(\lambda_s+\cdots+\lambda_d-\ell_r)\cr &\ge
\ell_0+\cdots+\ell_r-\ell_r\ge \ell_0+\cdots+\ell_i\cr}$$ et l'hypoth\`ese
de r\'ecurrence entra\^ine alors
$\sigma_{\bar\lambda'}\cdot\sigma_{\ell_0}\cdots\sigma_{\ell_{r-1}}\ne 0$,
ce qui montre le lemme.\cqfd  \bigskip

{\bf 5. Le r\'esultat de Hansen}

\ind Hansen a obtenu dans [H] un \theo\ de connexit\'e pour les
vari\'et\'es de drapeaux, qui s'\'enonce comme suit dans le cas des
grassmanniennes: {\it soient $X$ une vari\'et\'e irr\'eductible compl\`ete,
$\Delta$ la diagonale de $G(d,\P^n)\times G(d,\P^n) $ et $f:X\ra
G(d,\P^n)\times G(d,\P^n)$ un morphisme. Si $\dim f(X)< n$, alors
$f^{-1}(\Delta)$ est connexe.}

\ind Il montre aussi par la construction suivante que sa borne est la
meilleure possible.

(5.1) {\bf L'exemple de Hansen-Harris}. Dans $\P^n$, on consid\`ere une
conique lisse $C$ et un hyperplan $H$ la rencontrant transversalement. La
vari\'et\'e $X=\{\ (p,u)\in C\times G(d,\P^n)\mid p\in\Lambda_u \}$ est
projective irr\'eductible (et m\^eme lisse); notons $f:X\ra G(d,\P^n)$ la
seconde projection et $Y=G(d,H)$. Lorsque $d<n-1$, $f^{-1}(Y)$ a deux
composantes connexes. On a d'autre part $\codim f(X)=n-d-1$ et $\codim
(Y)=d+1$; plus pr\'ecis\'ement (\cf\ [Fu], ex. 14.7.6) $$
[f(X)]=2\sigma_{n-d-1}\qquad{\rm et}\qquad [Y]=\sigma_{ 1,\ldots,1}\ .$$
\ind On remarquera que $[f(X)]\cdot [Y]\cdot
(\sigma_{1,\ldots,1}+\sigma_{n-d})=0$.

(5.2) Dans la m\^eme veine, soient $d$ et $r$  des entiers tels que
$d/2+1\le r\le d+1$ et $d>0$. Posons $n=d+2r$ et consid\'erons une
quadrique lisse $Q$ dans $\P^n$. La vari\'et\'e $X=\{\ u\in G(d,\P^n)\mid
\Lambda_u\subset Q \}$ est irr\'eductible lisse ([H], th. 22.13). Soit $L$
un sous-espace lin\'eaire de $\P^n$ de dimension $(2r-1)$, transverse \`a
$Q$; soit $Y$ la vari\'et\'e de Schubert $\{ u\in G(d,\P^n)\mid
\dim(\Lambda_u\cap L)\ge r-1 \}$. Pour tout
 $u\in X\cap Y$, l'espace lin\'eaire $\Lambda_u\cap L$ est de dimension
$\ge r-1$ et est contenu dans la quadrique $L\cap Q$, lisse de dimension
$2r-2$. Il fait donc partie de l'une des deux familles de $(r-1)$\tx plans
contenus dans $L\cap Q$ (\loc ). On en d\'eduit que $X\cap Y$ n'est pas
connexe, alors que $\dim (X)+\dim (Y)>\dim G(d,\P^n)$. On notera que ([Fu],
ex. 14.7.15)  $$\nospacedmath [X]=2^{d+1}\sigma_{d+1,d,\ldots,1}\qquad{\rm
et}\qquad [Y]=\sigma_{\underbrace{\scriptstyle r,\ldots,r}_{ r\;\rm fois}}\
.$$ \ind Il ressort du lemme 4.1 que de nouveau, on a $[X]\cdot [Y]\cdot
(\sigma_{1,\ldots,1}+\sigma_{n-d})=0$.

\ind Notre but est d'\'etendre le r\'esultat de Hansen en tenant compte des
propri\'et\'es num\'eriques de la sous-vari\'et\'e $f(X)$ de
$G(d,\P^n)\times G(d,\P^n)$.

\bigskip {\bf 6. Un \theo\ de Bertini}

\ind Il faut \^etre prudent en ce qui concerne les \theo s du type Bertini
dans les grassmanniennes: dans  (5.2), l'intersection $X\cap \,\!^\gamma Y$
est non connexe pour $\gamma$ g\'en\'eral dans $\PGL(n+1,k)$, bien que
$\dim (X)+\dim(Y)>\dim G(d,\P^n)$. Dans [K], Kleiman construit, lorsque $k$
est de caract\'eristique non nulle, un exemple de surface dans $G(1,\P^3)$
dont l'intersection avec la vari\'et\'e de Schubert $G(1,H)$, o\`u $H$ est
un hyperplan {\it g\'en\'eral} de $\P^3$, est non r\'eduite.

\ind Le \theo\ suivant est notre r\'esultat cl\'e. Pour tout $p\in\P^n$, on
note $\Sigma_p$ la vari\'et\'e $\{ z\in G(1,\P^n)\mid p\in\Lambda_z \}$. Le
plongement  de Pl\"ucker induit un isomorphisme de $\Sigma_p$ sur un espace
projectif de dimension $n-1$. \smallskip {\pc TH\'EOR\`EME} 6.1.-- {\it
Soient $X$ une vari\'et\'e irr\'eductible,  $f:X\ra G(1,\P^n)$ un morphisme
dominant et $p$  un point g\'en\'eral de $\P^n$.

\ind {\rm a)}  $f^{-1}(\Sigma_p)$ est irr\'eductible;

\ind {\rm b)} $(k=\C)$ si $X$ est localement irr\'eductible, $\pi_1\bigl(
f^{-1}(\Sigma_p)\bigr) \dra \pi_1(X)$.} \medskip \ind Lorsque $X$ est
compl\`ete, le \theo\ de Hansen entra\^ine que   $f^{-1}(\Sigma_p)$ est
connexe pour tout $p$; il n'est pas difficile d'en d\'eduire a) dans ce
cas. On notera que l'hypoth\`ese que $f$ est dominant est indispensable: si
l'on prend $d=n-2$ dans l'exemple (5.1) et que l'on dualise, l'image de $f$
est un diviseur, mais $f^{-1}(\Sigma_p)$ a deux composantes connexes pour
$p$ g\'en\'eral.\medskip

{\bf D\'emonstration du \theo }. A ma grande surprise, je n'ai pas r\'eussi
\`a donner une d\'emonstration de ce r\'esultat bas\'ee sur le \theo\ de
Bertini usuel. Je vais donner deux d\'emonstrations: l'une sera valable en
toute caract\'eristique et d\'emontrera a), l'autre supposera $k=\C$ et
d\'emontrera a) et b). Posons $G=G(1,\P^n)$.

\ind Supposons donc tout d'abord $k$ de caract\'eristique quelconque. La
d\'emonstration suit celle du th. 6.3 de [J]. On rappelle qu'une extension
de corps $K\i K'$ est dite {\it primaire} si la fermeture alg\'ebrique de
$K$ dans $K'$ est radicielle. Prenons des coordonnees homog\`enes
$(x_0,\ldots,x_n)$ dans $\P^n$ et notons $\Lambda$ le sous-espace
lin\'eaire d\'efini par $x_0=x_1=0$. L'ouvert  $$G^0=\{\ u\in G\mid
\Lambda_u\cap\Lambda=\vide\ \}$$  de $G$ est affine; on peut le
param\'etrer en associant \`a $(a_2,\ldots,a_n,b_1,\ldots,b_n)\in k^{2n-2}$
la droite joignant les points $(1,0, a_2,\ldots,a_n)$ et
$(0,1,b_2,\ldots,b_n)$. Consid\'erons la correspondance d'incidence $I=\{
(p,u)\in  \A^{n+1}\times G^0\mid u\in\Sigma_{[p]} \}$; la condition
$u\in\Sigma_{[p]}$ est \'equivalente \`a $p\in\Lambda_u$, de sorte que $I$
est d\'efinie par les \'equations $p_i=p_0a_i+p_1b_i$, $i=2,\ldots,n$. En
particulier, l'application $$(p_0,p_1,a_2,\ldots,a_n,b_2,\ldots,b_n)\mapsto
(p_0,p_1,p_0a_2+p_1b_2,\ldots,p_0a_n+p_1b_n,a_2,\ldots,b_2,\ldots)$$
d\'efinit un isomorphisme de $G^0$\tx sch\'emas entre $\A ^2\times G^0$ et
$I$.

\ind Pour montrer le \theo , on peut remplacer $G$ par $G^0$. On note
$Z=X\times_{G^0}I$;  le morphisme $Z\ra X$ est un fibr\'e vectoriel trivial
de rang $2$. En particulier, $Z$ est int\`egre. Il s'agit de montrer que
les fibres du morphisme $h:Z\ra I\buildrel{pr_1}\over{\lra }\A^{n+1}$ sont
presque toutes irr\'eductibles. Le premier pas de l'argument est classique:
la fibre g\'en\'erique $F$ de $h$ est int\`egre, et il suffit de montrer
qu'elle est g\'eom\'etriquement irr\'eductible ([J], 4.10); pour cela,
montrons que le corps des fractions $K(Z)$ est extension primaire de
$K=k(p_0,\ldots,p_n)$ ([J], 4.3). Par hypoth\`ese, le morphisme $Z\ra I$
est dominant; on a une suite d'extensions $$\diagram{k(p_0,p_1)&\subset &
k(p_0,\ldots,p_n)&\subset &k(p_0,\ldots,p_n,a_2,\ldots,a_n) &\subset
&k(p_0,p_1,a_2,\ldots,b_2,\ldots)\ \
  \cr ||  &&||&&||&&||\ \ \quad\cr
K_0&\subset&K&\subset&L&\subset&\qquad\qquad K(I)\ \ \subset\ \  K(Z)\ ,\cr
}$$ o\`u $K(Z)$ est une extension pure de $K(X)$ de base $(p_0,p_1)$, et
$a_i$, $b_i\in K(X)$.

\ind Notant $\tilde K$ la fermeture s\'eparable de $K$ dans $K(Z)$, il
s'agit de montrer que $K=\tilde K$. Soit $\tilde L$ la fermeture
s\'eparable de $L$ dans $K(Z)$, on a $\tilde K\i\tilde L$ et $\tilde K$
(resp. $\tilde L$) est de type fini alg\'ebrique, donc fini sur $K$ (resp.
$L$). Pour tout $c\in k$, on d\'efinit un $K(X)$\tx automorphisme
$\sigma_c$ de $K(Z)$ par les relations $\sigma_c(p_0)=p_0+c$ et
$\sigma_c(p_1)=p_1$. On a $\sigma_c(p_i)=p_i+ca_i$ pour $i=2,\ldots,n$, de
sorte que $\sigma_c(L)\i L$ et $\sigma_c(\tilde L)\i \tilde L$. Le corps
$M_c=L(\sigma_c\tilde K)$ est une extension s\'eparable finie de $L$,
contenue dans $\tilde L$; comme $\tilde L$ est s\'eparable fini sur $L$, et
que $k$ est infini, il existe $c\ne c'$ tels que $M=M_c=M_{c'}$. On pose
$\theta_i=\sigma_c(p_i)=p_i+ca_i $ et
$\theta'_i=\sigma_{c'}(p_i)=p_i+c'a_i$, pour $i=2,\ldots,n$, de sorte que
$L=K_0(\theta_2,\ldots,\theta_n,\theta'_2,\ldots,\theta'_n)$. Soit $\xi$ un
\'el\'ement primitif de l'extension $\tilde K$ de $K$; posons
$a=\sigma_c(\xi)$ et $a'=\sigma_{c'}(\xi)$, de sorte que $M=L(a)=L(a')$.

\ind Puisque $X$ est g\'eom\'etriquement irr\'eductible, l'extension $k\i
K(X)$ est primaire; par [J], 3.13.3, il en est de m\^eme de l'extension
$K_0\i K(Z)$, et {\it a fortiori} de l'extension $K_0\i
K_0(a,\theta_2,\ldots,\theta_n)=\sigma_c\tilde K$. Puisque le morphisme
$I\ra\A^{2n}$ d\'efini par\break  $(p,u)\mapsto
(p_0,\ldots,p_n,a_2,\ldots,a_n)$ est dominant, le degr\'e de transcendance
de $L$ (donc aussi celui de $L(a)$) sur $K_0$ vaut $2n-2$; comme celui de
$\sigma_c\tilde K$ est $n-1$, la famille $\{ \theta'_2,\ldots,\theta'_n\}$
est alg\'ebriquement ind\'ependante sur $ \sigma_c\tilde K$. Par \loc ,
l'extension $K_0(\theta'_2,\ldots,\theta'_n)\i L(a)=M$ est primaire. Mais
$a'=\sigma_{c'}(\xi)$ appartient \`a $M$ et est alg\'ebrique s\'eparable
sur\break  $K_0(\theta'_2,\ldots,\theta'_n)=\sigma_{c'}K$. Par suite,
$a'=\sigma_{c'}(\xi)$ est dans $\sigma_{c'}K$, de sorte que $\xi$ est dans
$K$. On a donc $\tilde K=K$, d'o\`u le \theo .

\ind Supposons maintenant $k=\C$; on suit la d\'emonstration du th. 1.1 de
[FL1]. Par \og droite contenue dans $G$\fg , on entend toute
sous-vari\'et\'e de $G$ du type $$\ell_{p,P}=\{\ x\in G(1,\P^n)\mid
p\in\Lambda_x\subset P\ \}\ ,$$ o\`u $p$ est un point d'un $2$\tx plan $P$
contenu dans $\P^n$ (ce sont en fait toutes les droites contenues dans
l'image du plongement de Pl\"ucker de $G$). Elles sont param\'etr\'ees par
une vari\'et\'e de drapeaux ${\cal L}$ irr\'eductible lisse de dimension
$3n-4$.

\ind Comme dans la d\'emonstration du lemme 1.4 de [De], quitte \`a
r\'etr\'ecir $X$, on peut supposer que $f$ fait de $X$ un espace fibr\'e
topologiquement localement trivial au-dessus de $f(X)$, et que ce dernier
est le compl\'ementaire $G^1$ d'une hypersurface $B$ de $G$. Soient $u_0$
un point g\'en\'eral de $G^1$ et ${\cal L}_0$ la sous-vari\'et\'e de ${\cal
L}$ form\'ee des droites contenues dans $G$ et passant par
 $u_0$. Le sous-ensemble de ${\cal L}_0$ qui consiste en les droites
transverses \`a $B$ contient un ouvert de Zariski dense ${\cal L}_0^1$
([K]).
 Consid\'erons la correspondance d'incidence $ J=\{ (u,\ell)\in
G^1\times{\cal L}_0\mid u\in\ell \}$, et posons $Z=X\times_{G^1}J=\{
(x,\ell)\in X\times{\cal L}_0\mid f(x)\in\ell \}$.

\ind  L'image de la premi\`ere projection $J\ra G^1$ est $pr_1(J)=\{ u\in
G^1\mid \Lambda_u\cap\Lambda_{u_0}\ne\vide\}$, et r\'ealise $J$ comme
l'\'eclatement de $u_0$ dans $pr_1(J)$. On en d\'eduit que la premi\`ere
projection $Z\ra X$  r\'ealise $Z$ comme l'\'eclatement de $f^{-1}(u_0)$
dans $X_0=\{ x\in X\mid \Lambda_{f(x)}\cap\Lambda_{u_0}\ne\vide\}$.

\ind Posons $T=\{ (x,p)\in X\times\P^n\mid p\in\Lambda_{f(x)}\}$; la
premi\`ere projection $pr_1:T\ra X$ est un $\P^1$\tx fibr\'e, de sorte que
$T$ est irr\'eductible. La seconde projection $pr_2:T\ra\P^n$ est dominante
puisque $f$ l'est; puisque $u_0$ est g\'en\'eral, le \theo\  de Bertini
usuel entra\^ine que $pr_2^{-1}(\Lambda_{u_0})$ est irr\'eductible, donc
aussi $X_0=pr_1\bigl( pr_2^{-1}(\Lambda_{u_0})\bigr)$. Par suite, $Z$ est
irr\'eductible.

\ind La seconde projection $Z\ra J$ est un rev\^etement topologique; il en
est de m\^eme de la projection $J\ra {\cal L}_0$ au-dessus de l'ouvert
${\cal L}_0^1$ (ses fibres sont des sph\`eres priv\'ees de $\deg (B)$
points), donc aussi de la compos\'ee $h:Z\ra{\cal L}_0$ au-dessus de ${\cal
L}_0^1$. Or $h^{-1}({\cal L}_0^1)$, ouvert dans $Z$, est irr\'eductible, et
$h$ admet comme section $\ell\mapsto (x_0,\ell)$, pour tout point fix\'e
$x_0$ de $f^{-1}(u_0)$. Comme les fibres d'une fibration localement
triviale entre espaces connexes, qui admet une section, sont toutes
connexes, on en d\'eduit que pour tout $\ell\in {\cal L}_0^1$,
$f^{-1}(\ell)$ est connexe. Pour tout $\Sigma_p$ contenant une droite
transverse \`a $B$, $f^{-1}(\Sigma_p)$ est connexe; \'etant lisse, il est
irr\'eductible, ce qui d\'emontre a). Pour d\'emontrer b), il suffit
d'appliquer a) \`a la compos\'ee $f\circ\pi$, o\`u $\pi:\tilde X\ra X$ est
le rev\^etement universel de $X$.\cqfd \smallskip

\ind On en d\'eduit un \theo\ de Bertini pour les grassmanniennes. \bigskip
{\pc TH\'EOR\`EME} 6.2.-- {\it Soient $X$ une vari\'et\'e irr\'eductible,
$f:X\ra G(d,\P^n)$ un morphisme, et $l$ un entier tel que $d<l\le n$. On
suppose que pour $M\in G(l-1,\P^n)$ g\'en\'eral, $f(X)$ rencontre $G(d,M)$.
Pour $L\in G(l,\P^n)$ g\'en\'eral,

\ind {\rm a)}  $f^{-1}\bigl( G(d,L)\bigr)$ est irr\'eductible;

\ind {\rm b)} $(k=\C)$ si $X$ est localement irr\'eductible, $\pi_1\bigl(
f^{-1}\bigl( G(d,L)\bigr)\bigr) \dra \pi_1(X)$.}

\ind Le \theo\ de Bertini usuel correspond au cas $d=0$.

{\bf D\'emonstration}. Il suffit de traiter le cas o\`u $L$ est un
hyperplan de $\P^n$, o\`u l'on a $d\le n-2$. Posons  $$Z=\{\ (x,u)\in
X\times G(n-2,\P^n)\mid \Lambda_{f(x)}\subset \Lambda_u\ \}\ .$$ \ind La
premi\`ere projection $p:Z\ra X$ est un fibr\'e en grassmanniennes, de
sorte que $Z$ est irr\'eductible, et localement irr\'eductible si $X$
l'est. Dans le cadre topologique, on a aussi $\pi_1(Z) \dra \pi_1(X)$. La
seconde projection $q:Z\ra G(n-2,\P^n)$ est dominante par hypoth\`ese, et,
pour tout hyperplan $L$ de $\P^n$, on a $p\bigl( q^{-1}\bigl(
G(n-2,L)\bigr)\bigr)=f^{-1}\bigl( G(d,L)\bigr)$. On est donc ramen\'e au
cas $d=n-2$. En dua\-li\-sant, on obtient un morphisme $f:X\ra G(1,\P^n)$
dominant, auquel il suffit d'appliquer le \theo\ 6.1.\cqfd \medskip \ind
Avec des hypoth\`eses de propret\'e, on passe alors comme d'habitude \`a
des \'enonc\'es de connexit\'e.  \medskip
 {\pc TH\'EOR\`EME} 6.3.-- {\it Soient $X$ une vari\'et\'e irr\'eductible,
$f:X\ra G(d,\P^n)$ un morphisme, $l$ un entier tel que $d<l\le n$, et $L$
un sous-espace lin\'eaire de $\P^n$ de dimension $l$. On suppose que pour
$M\in G(l-1,\P^n)$ g\'en\'eral, $f(X)$ rencontre $G(d,M)$.

\ind {\rm a)}  Si $f$ est propre au-dessus d'un ouvert $V$ de $G(d,\P^n)$
et que $G(d,L)$ est contenu dans $V$, alors $f^{-1}\bigl( G(d,L)\bigr)$ est
connexe;

\ind {\rm b)} $(k=\C)$ si $X$ est localement irr\'eductible, pour tout
voisinage ouvert $U$ de $G(d,L)$ dans $G(d,\P^n)$, on a $\pi_1\bigl(
f^{-1}(U)\bigr)\dra \pi_1(X)$.}

{\bf D\'emonstration}. On suit [FL1], th. 2.1; posons $Z=\{ (x,u)\in X\times
G(d,L) \mid \Lambda_{f(x)}\subset \Lambda_u \}$. La premi\`ere projection
r\'ealise $Z$ comme un fibr\'e en grassmanniennes au-dessus de $X$, de
sorte que $Z$ est irr\'eductible. La projection $q:Z\ra G(d,L)$ est propre,
donc admet une factorisation de Stein $Z\buildrel q'\over{\ra}G'\buildrel
\rho\over{\ra} G(d,L)$, o\`u $\rho$ est fini et o\`u les fibres de $q'$
sont connexes. Par th. 6.2, $\rho$ est birationnel surjectif; comme
$G(d,L)$ est normal, $\rho$ est bijectif, de sorte que les fibres de $q$
sont connexes, ce qui montre a).

\ind Tout voisinage ouvert $U$ de $G(d,L)$ contient des $G(d,L')$ auxquels
on peut appliquer le th. 6.2.b), d'o\`u b).\cqfd \medskip {\it Remarque}
6.4.-- On aura besoin de la version plus fine de b) suivante: {\it pour
tout $x\in f^{-1}\bigl( G(d,L)\bigr)$, l'homomorphisme $\pi_1\bigl(
f^{-1}(U),x\bigr)\ra \pi_1(X,x)$ est surjectif}. Cela se d\'e\-mon\-tre
comme dans [FL1], remark 2.2.
 \bigskip {\bf 7. Th\'eor\`emes de connexit\'e}

\ind On note $\Delta$ la diagonale de $G(d,\P^n)\times G(d,\P^n) $ et, pour
tout
 $\gamma\in \PGL (n+1,k)$, $^\gamma\!\Delta$ l'image de $\Delta$ par
l'automorphisme $(x,y)\mapsto (x,\gamma y)$ de $G(d,\P^n)\times G(d,\P^n)$.
\smallskip {\pc TH\'EOR\`EME} 7.1.-- {\it Soient $X$ une vari\'et\'e
irr\'eductible et $f:X\ra G(d,\P^n)\times G(d,\P^n)$ un morphisme. On
suppose qu'il existe des partitions $\lambda$ et $\mu$ v\'erifiant
$\lambda<\bar\mu$ ou $\lambda^*<\bar\mu^*$, telles que $
 [\overline {f(X)}]\cdot p_1^*\sigma_{\bar\lambda}\cdot
p_2^*\sigma_{\bar\mu}\ne 0$.

\ind {\rm a)} Pour $\gamma$ g\'en\'eral dans $\PGL (n+1,k)$,
$f^{-1}(\,\!^\gamma\!\Delta)$ est irr\'eductible;

\ind {\rm b)} si $X$ est compl\`ete, $f^{-1}(\Delta)$ est connexe;

\ind {\rm c)} $(k=\C)$  si $X$ est localement irr\'eductible et compl\`ete,
$\pi_1\bigl( f^{-1}(\Delta)\bigr)\dra \pi_1(X)$.}

{\it Remarques} 7.2.-- 1) La condition $\lambda<\bar\mu$ \'equivaut aux
in\'egalit\'es $\lambda_i+\mu_{d-i}< n-d$ pour tout $i=0,\ldots,d$. La
condition $\lambda^*<\bar\mu^*$ \'equivaut aux relations
$\lambda_d=\mu_d=0$ et $\lambda_i+\mu_{d-i-1}\le n-d$ pour tout
$i=0,\ldots,d-1$.

\ind 2) Le \theo\ entra\^ine le r\'esultat de Hansen (\cf\ \S 5); en effet,
supposons $\codim f(X)<n$ et soient $\lambda$ et $\mu$ deux partitions
v\'erifiant $[f(X)]\cdot p_1^*\sigma_{\bar\lambda}\cdot
p_2^*\sigma_{\bar\mu}\ne 0$ et $|\lambda|+|\mu|=\codim f(X)$. Si
$\lambda_\alpha\ge\bar\mu_\alpha$, on a $$\eqalign{|\lambda|+|\mu| &\ge
(\alpha+1)\lambda_\alpha+(d+1-\alpha)(n-d-\lambda_\alpha)\cr
&=\lambda_\alpha(2\alpha-d)+(d+1-\alpha)(n-d)\ .\cr}$$ \ind Quitte \`a
\'echanger $\lambda$ et $\mu$, on peut supposer $\alpha\le d/2$. Si
$\lambda_\alpha<n-d$, on a $$\eqalign{n>\codim f(X)= |\lambda|+|\mu|&\ge
(n-d-1)(2\alpha-d)+(d+1-\alpha)(n-d)\cr
&=(n-d)(\alpha+1)-2\alpha+d=n+(n-d-2)\alpha\ ,\cr}$$ de sorte que $n-d\le
1$, $\lambda_\alpha=0$ et $\mu_{d-\alpha}=n-d$.

\ind On a donc dans tous les cas $\lambda_0=n-d$ ou $\mu_0=n-d$. Si
l'in\'egalit\'e $\lambda^*<\bar\mu^*$ est aussi viol\'ee, on obtient de
m\^eme $\lambda_0^*=d+1$ ou $\mu_0^*=d+1$, d'o\`u la contradiction
$$|\lambda|+|\mu| \ge n-d +(d+1)-1=n \ .$$ \medskip   {\bf D\'emonstration
du \theo }.  On utilise de nouveau l'astuce de Deligne pour se ramener au
\theo\ de Bertini 6.2. On consid\`ere des coordonn\'ees homog\`enes
$[x^{(1)},x^{(2)}]$ sur $\P^{2n+1}$, o\`u $x^{(1)}$ et $x^{(2)}$ sont des
$(n+1)$\tx uplets d'\'el\'ements de $k$. Soit $\tilde\gamma$ un \'el\'ement
de $\GL (n+1,k)$ et $\gamma$ son image dans $\PGL (n+1,k)$; on note $L_1$
(resp.~$L_2$) (resp.~$^{\tilde\gamma}L$) le sous-espace lin\'eaire de
dimension $n$ de $\P^{2n+1}$ d\'efini par $x^{(1)}= 0$  (resp.~$x^{(2)}=0$)
(resp.~$x^{(1)}=\tilde\gamma(x^{(2)})$). Soit $V$ l'ouvert
$G(d,\P^{2n+1})\moins\Sigma_{L_1}\moins\Sigma_{L_2}$; pour $i=1$, $2$,
soient $\rho_i:\P^{2n+1}\moins L_i\ra\P^n$ le morphisme
$[x^{(1)},x^{(2)}]\mapsto [x^{(i)}]$ et $\phi_i:V\ra  G(d,\P^n)$ le
morphisme induit. Le morphisme  $\phi=(\phi_1,\phi_2):V\ra G(d,\P^n)\times
G(d,\P^n)$ est un\break $\GL(d+1,k)$\tx fibr\'e, et induit un isomorphisme
de $G(d,\!\,^{\tilde\gamma}L)$ sur
 $^\gamma\!\Delta$. Si on note\break $X'=X\times_{G(d,{\bf P}^n)\times
G(d,{\bf P}^n)}V$ et $f':X'\ra V$ la projection, on a $f'^{-1}\bigl(
G(d,\!\,^{\tilde\gamma}L)\bigr)\isom f^{-1}(\!\,^\gamma\!\Delta)$.

\ind Soit $M$ un sous-espace lin\'eaire de $\P^{2n+1}$  g\'en\'eral de
dimension $n-1$, de sorte que, pour $i=1$, $2$, l'espace $M\cap L_i$ est
vide et que $\rho_i$ induit un isomorphisme de $M$ sur un hyperplan $H_i$ de
$\P^n$. On note $\gamma_i:G(d,M)\ra G(d,H_i)$ l'isomorphisme induit. Le
morphisme $\phi$ induit un isomorphisme de $G(d,M)$ sur $$\Omega_M=\{\
(u,\gamma_2\bigl( \gamma_1^{-1}(u)\bigr)\mid u\in G(d,H_1)\ \}$$ et
$f'^{-1}\bigl( G(d,M)\bigr)\isom f^{-1}(\Omega_M)$. La d\'ecomposition de
K\"unneth de $A^*\bigl(G(d,\P^n)\times G(d,\P^n)\bigr)$ ([FMSS]) permet de
d\'eterminer la classe de la sous-vari\'et\'e $\Omega_M$; elle vaut $$\sum
p_1^*\sigma_\alpha\cdot p_2^*\sigma_\beta\ ,$$ la somme portant sur toutes
les partitions $\alpha$ et $\beta$ telles que $\alpha_i+\beta_{d-i}=n-d+1$
pour tout $i=0,\ldots,d$.

 \ind Quitte \`a dualiser,  il existe par hypoth\`ese des partitions
$\lambda$ et $\mu$ telles que $\lambda<\bar\mu$, v\'erifiant $
[\overline{f(X)}]\cdot p_1^*\sigma_{\bar\lambda}\cdot
p_2^*\sigma_{\bar\mu}\ne 0$. D\'efinissons une partition $\lambda'$ par
$\lambda'_i=\bar\mu_i-1$; on a $\lambda'\ge\lambda$, et le lemme 4.1
entra\^ine $ [\overline{f(X)}]\cdot p_1^*\sigma_{\bar\lambda'}\cdot
p_2^*\sigma_{\bar\mu}\ne 0$. Comme le produit $
p_1^*\sigma_{\bar\lambda'}\cdot p_2^*\sigma_{\bar\mu}$ intervient dans
l'expression de $[\Omega_M]$ donn\'ee ci-dessus, on obtient
$[\overline{f(X)}]\cdot [\Omega_M]\ne 0$. On a donc
$[\overline{f'(X')}]\cdot [G(d,M)]\ne 0$, d'o\`u a) (th. 6.2). Lorsque $X$
est compl\`ete, $f'$ est propre; on en d\'eduit b) (th. 6.3).

\ind Le point c) se d\'eduit du th. 6.3.b) et de la remarque 6.4, comme
dans [FL1], p. 40.\cqfd \medskip \ind On peut d\'erouler les corollaires
habituels. \medskip

{\pc COROLLAIRE} 7.3.-- {\it Soient $X$ et $Y$ des vari\'et\'es
irr\'eductibles
 et $f:X\ra G(d,\P^n)$ et $g:Y\ra G(d,\P^n)$ des morphismes. On suppose que
$[\overline{f(X)}]\cdot [\overline{g(Y)}]\cdot
(\sigma_{1,\ldots,1}+\sigma_{n-d})\ne 0$.

\ind {\rm a)} Pour $\gamma$ g\'en\'eral dans $\PGL (n+1,k)$,
$X\times_{G(d,{\bf P}^n)} \,\! ^\gamma Y$ est irr\'eductible;

\ind {\rm b)} si $X$ et $Y$ sont compl\`etes, $X\times_{G(d,{\bf P}^n)} Y$
est connexe;

\ind {\rm c)} $(k=\C)$ si $X$ et $Y$ sont localement irr\'eductibles et
compl\`etes,\break $\pi_1(X\times_{G(d,{\bf P}^n)} Y)\dra \pi_1(X\times
Y)$.} \smallskip {\bf D\'emonstration}. Il suffit de remarquer que pour que
$[\overline{f(X)}]\cdot [\overline{g(Y)}]\cdot \sigma_{1,\ldots,1}$ soit
non nul, il faut et il suffit qu'il existe des partitions $\lambda$ et
$\mu$ v\'erifiant $\lambda<\bar\mu$, telles que
$[\overline{f(X)}]_\lambda\ne 0$ et $[\overline{g(Y)}]_\mu\ne 0$ (lemme
4.1). On applique ensuite le \theo\ 7.1 au morphisme $(f,g):X\times Y\ra
G(d,\P^n)\times G(d,\P^n)$.\cqfd

\ind Lorsque $k=\C$, que $X$ et $Y$ sont compl\`etes, que $f$ est {\it
surjective} et que $g$ est non constante, on peut montrer en utilisant
[PS], prop. 1, que les conclusions du \theo\ subsistent, bien que
l'hypoth\`ese $[g(Y)]\cdot (\sigma_{1,\ldots,1}+\sigma_{n-d})\ne 0$ ne soit
pas toujours v\'erifi\'ee.

 \bigskip {\pc COROLLAIRE} 7.4.-- {\it Soit $X$ une sous-vari\'et\'e
irr\'eductible de $G(d,\P^n)$ telle que\break $[X]\cdot [X]\cdot
(\sigma_{1,\ldots,1}+\sigma_{n-d})\ne 0$.

\ind {\rm a)} $\pi_1^{\rm alg}(X)=0$;

\ind {\rm b)} $(k=\C)$  $X$ est simplement connexe.} \bigskip

\ind Rappelons (\S 2) qu'une sous-vari\'et\'e  $Z$ de $G(d,\P^n)$ est {\it
encombrante} si elle rencontre toute sous-vari\'et\'e de $G(d,\P^n)$ de
dimension $\ge\codim (Z)$. De fa\c con \'equivalente, on demande que pour
toute partition $\lambda$ telle que $|\lambda|=\codim(Z)$, on ait
$[Z]_\lambda\ne 0$.

\ind Soit $Q$ le fibr\'e quotient universel sur $G(d,\P^n)$; il d\'ecoule de
[FL2] que toute une sous-vari\'et\'e $Z$ de $G(d,\P^n)$ telle que la
restriction de $Q$ \`a $Z$ soit ample, est encombrante (ainsi par
cons\'equent que toute sous-vari\'et\'e de $Z$).  C'est le cas par exemple
pour l'image $Z$ de l'application de Gauss d'une sous-vari\'et\'e lisse
d'une vari\'et\'e ab\'elienne {\it simple} (\cf\ [D2]). En consid\'erant
des intersections compl\`etes de codimension $2$ dans une vari\'et\'e
ab\'elienne simple de dimension $n+1$, on obtient des exemples de
sous-vari\'et\'es encombrantes de dimension $n-1$ dans $G(1,\P^n)$ dont le
groupe fondamental est isomorphe \`a $\Z^{2n+2}$.\bigskip

{\pc COROLLAIRE} 7.5.-- {\it Soient $X$ une vari\'et\'e irr\'eductible
compl\`ete et $f:X\ra G(d,\P^n)$ un morphisme dont l'image est une
sous-vari\'et\'e encombrante. Pour toute sous-vari\'et\'e irr\'eductible
$Z$ de $G(d,\P^n)$ de dimension  $>\codim f(X)+d$, telle que $[Z]\cdot
\sigma_{1,\ldots,1}\ne 0$, $f^{-1}(Z)$ est connexe.}

\ind La conclusion du corollaire subsiste lorsque $\dim (Z)>\codim
f(X)+n-d-1$ et  $[Z]\cdot\sigma_{n-d}\ne 0$: il suffit de dualiser (\cf\
aussi cor. 8.3). D'autre part, rappelons que si l'on prend $d=n-2$ dans
l'exemple (5.1), on obtient un morphisme $f:X\ra G(1,\P^n)$ dont l'image
est un diviseur (ample, donc encombrant), tel que $f^{-1}(\Sigma_{1,1})$
ait deux composantes connexes.  \medskip

{\bf D\'emonstration de la proposition}. Il existe une partition $\lambda$
avec  $\lambda_0<n-d$ et $[f(X)]_\lambda\ne 0$. Si on pose
$\mu_i=\bar\lambda_i-1$, on a $\lambda<\bar\mu$ et
$|\mu|=|\bar\lambda|-d>\codim (Z)$. Il existe une partition $\mu'$ telle
que $\mu'\le\mu$ et $|\mu'|=\codim (Z)$; comme $Z$ est encombrante, on a
$[Z]_\mu\ne 0$, puis $[Z]_{\mu'}\ne 0$ (lemme 4.1), et la proposition
r\'esulte du cor. 7.3.\cqfd

 \bigskip {\bf 8. Images inverses des vari\'et\'es de Schubert}

\ind Soit $\mu=(\mu_0,\ldots,\mu_d)$ une partition; on pose $J(\mu)=\{
j\in\{ 0,\ldots,d\} \mid \mu_j>\mu_{j+1} \}$. Si $\mu_{d}<n-d$, on
d\'efinit pour tout $j\in J(\mu)$ une partition $\mu^{(j)}$ de la fa\c con
suivante: si $\mu_j<n-d$, on pose $\mu^{(j)}_i=\mu_j+1$ pour $i\le j$, et
$\mu^{(j)}_i=\mu_i$ pour $i> j$; si $\mu_j=n-d$, on pose $\mu^{(j)}_i=n-d$
pour $i\le j+1$, et $\mu^{(j)}_i=\mu_i$ pour $i> j+1$. \medskip {\pc
TH\'EOR\`EME} 8.1.-- {\it Soient $X$ une vari\'et\'e irr\'eductible,
$f:X\ra G(d,\P^n)$ un morphisme, et $\Sigma_{\bf L}$ une vari\'et\'e de
Schubert de dimension $>0$ et de classe $\sigma_\mu$. On suppose que pour
tout $j\in J(\mu)$, on a $[\overline{f(X)}]\cdot\sigma_{\mu^{(j)}}\ne 0$.

\ind{\rm a)} Si ${\bf L}$ est g\'en\'eral, $f^{-1}(\Sigma_{\bf L})$ est
irr\'eductible;

\ind{\rm b)} si $X$ est compl\`ete, $f^{-1}(\Sigma_{\bf L})$ est connexe;

\ind{\rm c)}  $(k=\C)$ si $X$ est localement irr\'eductible et compl\`ete,
 $\pi_1\bigl( f^{-1}(\Sigma_{\bf L})\bigr)\dra \pi_1(X)$.}

{\it Exemples} 8.2.-- 1) Prenons par exemple $\mu=(5,2,2,1)$; lorsque
$n-d>5$, on demande que le produit de  $[\overline{f(X)}]$ avec chacune des
classes de Schubert de type $(2,2,2,2)$, $(3,3,3,1)$ et $(6,2,2,1)$ soit
non nul (pour $n-d=5$, changer la derni\`ere en $(5,5,2,1)$). Ce \theo\ est
donc en g\'en\'eral meilleur que le th. 7.1. On notera que les hypoth\`eses
du \theo\ ne sont pas invariantes par dualit\'e: sur l'exemple, la
condition duale est, pour $d\ge 4$, que le produit de  $[\overline{f(X)}]$
avec chacune des classes de Schubert de type $(5,5)$, $(5,2,2,2)$ et
$(5,2,2,1,1)$ soit non nul (lorsque $d=3$, supprimer la derni\`ere).

\ind 2) Lorsque $n-d>\mu_0>\cdots >\mu_r>\mu_{r+1}=0$, la condition requise
est que  l'inter\-sec\-tion de  $[\overline{f(X)}]$ avec chacune des
classes de Schubert de type
$(\mu_i+1,\ldots,\mu_i+1,\mu_{i+1},\ldots,\mu_r)$ soit non nulle
($i=0,\ldots,r$); il suffit pour cela que l'intersection de
$[\overline{f(X)}]$ avec  chacune des classes de Schubert qui apparaissent
dans  $\sigma_\mu\cdot\sigma_1$ soit non nulle. En particulier, pour une
vari\'et\'e de Schubert sp\'eciale de classe $\sigma_m$, on obtient la
condition  $[\overline{f(X)}]\cdot\sigma_{m+1}\ne 0$. La condition duale
est $[\overline{f(X)}]\cdot\sigma_{m,m}\ne 0$. Plus g\'en\'eralement, pour
une vari\'et\'e de Schubert de classe
$\nospacedmath\sigma_{\underbrace{\scriptstyle m,\ldots,m}_{ r\;\rm
fois}}$, les deux conditions peuvent se regrouper en $$\nospacedmath
[\overline{f(X)}]\cdot (\sigma_{\underbrace{\scriptstyle m+1,\ldots,m+1}_{
r\;\rm fois}}+\sigma_{\underbrace{\scriptstyle m,\ldots,m}_{ r+1\;\rm
fois}})\ne 0\ .$$ \ind Lorsque $m=r=1$, on obtient
$[\overline{f(X)}]\cdot(\sigma_2+\sigma_{1,1})=[\overline{f(X)}]\cdot\sigma_1^2\ne
0$, \cad\ simplement $\dim f(X)\ge 2$, l'hypoth\`ese du \theo\ de Bertini
usuel.

 \medskip {\bf D\'emonstration du \theo }. Notons la vari\'et\'e  $F$ des
drapeaux $(\Lambda_0\subset\cdots\subset\Lambda_d\subset\P^n)$, o\`u
$\Lambda_i$ a dimension $i$, et $p_i$ la surjection $F\ra G(i,\P^n)$.
Soient $J$ un sous-ensemble de $\{ 0,\ldots,d\}$ et $\{ M^{(j)}\}_{j\in J}$
une famille croissante de sous-espaces lin\'eaires de $\P^n$. On a
$$p_d\Bigl( \bigcap_{j\in J}p_j^{-1}\bigl( G(j,M^{(j)})\bigr)\Bigr)=\{\
u\in G(d,\P^n)\mid \dim (\Lambda_u\cap M^{(j)})\ge j\quad{\rm pour\ tout}\
\ j\in J\ \}\ .$$ \ind Si la fonction $j\mapsto \dim (M^{(j)})-j$ est
croissante et \`a valeurs dans $\{ 0,\ldots,n-d\}$, cette sous-vari\'et\'e
de $G(d,\P^n)$ est la  vari\'et\'e de Schubert associ\'ee \`a tout drapeau
$(N^{(0)}\subset\cdots\subset N^{(d)})$ o\`u $N^{(i)}=M^{(i)}$ si $i\in J$,
o\`u $N^{(i)}$ est un hyperplan quelconque de $N^{(i+1)}$ si $i\notin J$ et
o\`u  $N^{(d)}=\P^n$ si $d\notin J$. La partition correspondante $\lambda$
est d\'efinie par $\lambda_i=n-d+j-\dim  (M^{(j)})$, o\`u $j$ est le plus
petit \'el\'ement de $J$ sup\'erieur \`a $i$.

\ind Notons $(L^{(0)}\subset \cdots\subset L^{(d)} \subset\P^n)$ le drapeau
${\bf L}$;  la vari\'et\'e de Schubert $\Sigma_{\bf L}$ peut se d\'efinir
par les in\'egalit\'es $\dim (\Lambda_u\cap L^{(j)})\ge j$ en se
restreignant aux $j$ dans $J(\mu)$. En particulier, $$\Sigma_{\bf
L}=p_d\Bigl( \bigcap_{j\in J(\mu)}p_j^{-1}\bigl( G(j,L^{(j)})\bigr)\Bigr)\
.$$ \ind La vari\'et\'e $X'=X\times_{G(d,{\bf P}^n)}F$ est irr\'eductible
puisque $X$ l'est; notons $f'$ la projection $X'\ra F$. Pour $i\ge 0$,
notons $X'_i$ la sous-vari\'et\'e $\displaystyle\bigcap_{j\in J(\mu),\,
j\ge i}(p_jf')^{-1}\bigl( G(j,L^{(j)})\bigr)$ de $X'$, de sorte que
$X'_{d+1}=X'$, et que l'image de $X'_0$ par la projection $X'\ra X$ est
$f^{-1}(\Sigma_{\bf L})$.

\ind Supposons ${\bf L}$ g\'en\'eral et montrons par r\'ecurrence
descendante sur $i$ que $X'_i$ est irr\'eductible, ce qui prouvera a).
Supposons $X'_i$ irr\'eductible pour tout $i>j$ et montrons que $X'_j$
l'est aussi. Vue la d\'efinition de $X'_j$, il suffit de traiter le cas
o\`u $j\in J(\mu)$.

 \ind Supposons d'abord que $j$ ne soit pas le plus grand \'el\'ement de
$J(\mu)$, et soit $j'$ le plus petit \'el\'ement de $J(\mu)$ strictement
sup\'erieur \`a $j$. Le morphisme $p_j f'$ se restreint en un morphisme
$f_j:X'_{j'}\ra G(j,L^{(j')})$ et $X'_j=f_j^{-1}\bigl( G(j,L^{(j)})\bigr)$.
Le th. 6.2 montre que $X'_j$ est irr\'eductible sous les hypoth\`eses
suivantes: d'une part $\dim (L^{(j)})>j$, ce qui \'equivaut \`a
$\mu_j<n-d$, d'autre part $f'(X')$ rencontre  $p_j^{-1}\bigl(
G(j,M^{(j)})\bigr)$, o\`u $M^{(j)}$ est un sous-espace lin\'eaire de
$L^{(j')}$ g\'en\'eral de dimension $\dim (L^{(j)})-1$. Cette derni\`ere
condition \'equivaut \`a  $$f(X)\cap p_d\Bigl(p_j^{-1}\bigl(
G(j,M^{(j)})\cap \bigcap_{k\in J(\mu),\, k>j}p_k^{-1}\bigl(
G(k,L^{(k)})\bigr)\Bigr)\ne \vide\ ,$$
 soit encore, d'apr\`es la discussion pr\'ec\'edente, \`a l'hypoth\`ese
$[\overline{f(X)}]\cdot\sigma_{\mu^{(j)}}\ne 0$.

\ind Lorsque $\mu_j=n-d$, on consid\`ere l'application $f_{j'}:X'_{j'}\ra
G(j',L^{(j')})$; la vari\'et\'e $X'_j$ n'est autre que l'image inverse de
$\Omega_{L^{(j)}}=\{ u\in G(j',L^{(j')})\mid\Lambda_u\supset L^{(j)}\ \}$
par $f_{j'}$. La version duale du th. 6.2 entra\^ine que pour que $X'_j$
soit irr\'eductible, il suffit que $\dim (L^{(j)})<j'$ et que
$f_{j'}^{-1}(\Omega_{M^{(j)}})$ soit non vide pour tout sous-espace
lin\'eaire $M^{(j)}$ de $L^{(j')}$ de dimension $\dim (L^{(j)})+1$. De
nouveau, la discussion pr\'ec\'edente montre que cette condition est
\'equivalente \`a l'hypoth\`ese
$[\overline{f(X)}]\cdot\sigma_{\mu^{(j)}}\ne 0$.

\ind La d\'emonstration du pas de r\'ecurrence lorsque $j$ est le plus
grand \'el\'ement de $J(\mu)$ est identique: il suffit de remplacer
$X'_{j'}$ par $X'$ et $L^{(j')}$ par $\P^n$, puis $f_{(j')}$ par $p_d$ et
$j'$ par $d$. Ceci termine la d\'emonstration de a).

  \ind Lorsque $k=\C$ et que $X$ est localement irr\'eductible, ce qui
pr\'ec\`ede montre aussi que
 $\pi_1\bigl( f^{-1}(\Sigma_{\bf L})\bigr)\dra \pi_1(X)$ (toujours pour
${\bf L}$ g\'en\'eral). On en d\'eduit c) comme dans [FL1] rem. 2.2 et cor.
3.3.

\ind Montrons b). Soit $F'$ la vari\'et\'e de drapeaux contenant le drapeau
${\bf L}$; posons  $$Z=\{\ (x,{\bf L'})\in X\times F' \mid \dim
(\Lambda_{f(x)}\cap L'^{(i)})\ge i\ \ {\rm pour\ tout}\ \ i=0,\ldots,d\ \}\
.$$   \ind La premi\`ere projection r\'ealise $Z$ comme un fibr\'e
au-dessus de $X$, dont les fibres sont des vari\'et\'es de Schubert (donc
irr\'eductibles par [Fu], ex. 14.7.16) de m\^eme type dans la vari\'et\'e
de drapeaux $F'$, de sorte que $Z$ est irr\'eductible. Supposons $X$
compl\`ete; la projection $q:Z\ra F'$ est propre, donc admet une
factorisation de Stein $Z\buildrel q'\over{\ra}F''\buildrel \rho\over{\ra}
F'$, o\`u $\rho$ est fini et o\`u les fibres de $q'$ sont connexes. Par a),
$\rho$ est birationnel surjectif; comme $F'$ est normal, $\rho$ est
bijectif, de sorte que les fibres de $q$ sont connexes, ce qui montre
b).\cqfd \medskip \ind Pour toute partition $\mu$ et tout $j\in J(\mu)$, on
d\'efinit un entier $\delta(j)$ par: $$\eqalign{{\rm si}\quad
\mu_j<n-d\qquad\quad &\delta(j)=|\mu^{(j)}|-|\mu|-1=\sum_{i<
j}(\mu_j-\mu_i+1)\cr {\rm si}\quad \mu_j=n-d\qquad\quad
&\delta(j)=n-d-1-\mu_{j+1}\ ,\cr}$$ et on pose $\delta(\mu)=\sum_{j\in
J(\mu)}\max \bigl(\delta(j),0\bigr)$; c'est un \'el\'ement de $\{
0,\ldots,\max(d,n-d-1)\}$, qui peut prendre toutes les valeurs dans cet
ensemble. Lorsque $\mu_0<n-d$ et que la partition $\mu$ est strictement
d\'ecroissante, $\delta(\mu)=0$; lorsque  $\mu$ est constante,
$\delta(\mu)=\delta(\mu^*)=d$; lorsque $\mu=(n-d,0,\ldots)$,
$\delta(\mu)=\delta(\mu^*)=n-d-1$. \medskip  {\pc COROLLAIRE} 8.3.-- {\it
Soient $X$ une vari\'et\'e irr\'eductible compl\`ete et $f:X\ra G(d,\P^n)$
un morphisme dont l'image est une sous-vari\'et\'e encombrante. Si
$\Sigma_\mu$ est une sous-vari\'et\'e de Schubert de $G(d,\P^n)$, de
dimension $>\codim f(X)+\delta(\mu)$, $f^{-1}(\Sigma_\mu)$ est connexe.}

\ind Le lecteur remarquera que $\delta(\mu)$ et $\delta(\mu^*)$ peuvent
\^etre diff\'erents (c'est le cas par exemple lorsque $\mu=(3,3,1,1)$); on
aura donc parfois int\'er\^et \`a dualiser avant d'appliquer le corollaire.
D'autre part, rappelons que si l'on prend $d=n-2$ dans l'exemple (5.1), on
obtient un morphisme $f:X\ra G(n-2,\P^n)$ dont l'image est un diviseur
(ample, donc encombrant), tel que $f^{-1}(\Sigma_{1,1})$ ait deux
composantes connexes. La borne du corollaire est donc la meilleure possible.

\bigskip

\ind Nous terminerons avec un \'enonc\'e portant sur les images inverses
d'intersections de vari\'et\'es de Schubert sp\'eciales. Le r\'esultat
obtenu est en g\'en\'eral meilleur que le \theo\ g\'en\'eral 7.1.

 \medskip

{\pc TH\'EOR\`EME} 8.4.-- {\it Soient $X$ une vari\'et\'e irr\'eductible,
$f:X\ra G(d,\P^n)$ un morphisme et $L_0,\ldots,L_r$ des sous-espaces
lin\'eaires de $\P^n$. On note $\ell_i=  \codim \Sigma_{L_i}$ et on suppose
que $\ell_0\ge\ell_1\ge\cdots\ge\ell_r$. Soit $s$ le cardinal de l'ensemble
$\{ i\mid\ell_i=n-d\}$. On suppose que $[\overline{f(X)}]\cdot
\sigma_{n-d}^s\cdot\sigma_{\ell_s+1}\cdot\sigma_{\ell_{s+1}}
\cdots\sigma_{\ell_r}\ne 0$ si $s\le r$ et que $[\overline{f(X)}]\cdot
\sigma_{n-d}^{s+1}\ne 0$ si $s>0$.

\ind\ind{\rm a)} Si $L_0,\ldots,L_r$ sont g\'en\'eraux,
$f^{-1}(\Sigma_{L_0}\cap\cdots\cap \Sigma_{L_r})$ est irr\'eductible;

\ind\ind{\rm b)} si $X$ est  compl\`ete, $f^{-1}(\Sigma_{L_0}\cap\cdots\cap
\Sigma_{L_r})$ est connexe;

\ind\ind{\rm c)}  $(k=\C)$ si $X$ est localement irr\'eductible et
compl\`ete,\break
 $\pi_1\bigl( f^{-1}(\Sigma_{L_0}\cap\cdots\cap \Sigma_{L_r})\bigr)\dra
\pi_1(X)$.} \medskip {\it Remarques} 8.5.-- 1) Compte tenu du lemme 4.2,
l'hypoth\`ese du \theo\ est \'equivalente \`a la propri\'et\'e suivante: il
existe une partition $\lambda$ avec $[\overline{f(X)}]_\lambda\ne0$ telle
que, pour tout $i=0,\ldots,r$, on ait
$\ell_0+\cdots+\ell_i\le\bar\lambda_0+\cdots+\bar\lambda_i$, avec
in\'egalit\'e stricte pour $i\ge s$. Lorsque $s>0$, il faut ajouter \`a
cela la condition qu'il existe une partition $\lambda'$ avec
$\lambda'_{d-s}=0$ et $[\overline{f(X)}]_{\lambda'}\ne0$.

\ind 2) Posons $c=\codim f(X)$ et supposons $c\le n-d$; si
$[\overline{f(X)}]_{(c)}\ne 0$, l'hypoth\`ese du \theo\ est satisfaite
d\`es que  $\dim f(X)>\sum_{i=0}^r\ell_i$ (lemme 4.2); lorsque $s>0$, il
faut aussi supposer  $\dim f(X)\ge(s+1)(n-d)$. On g\'en\'eralise ainsi le
r\'esultat de [HS], sauf dans le cas o\`u $s>0$ et
$\sum_{i=s}^r\ell_i<n-d$. On notera par ailleurs que la d\'emonstration de
\loc\ de l'irr\'eductibilit\'e d'une intersection de vari\'et\'es de
Schubert sp\'eciales g\'en\'erales n'est valable qu'en caract\'eristique
nulle.

{\bf D\'emonstration du \theo }. {\it Supposons d'abord $s=0$};
d\'efinissons $$Z=\{\ (x,v_0,\ldots,v_r)\in X\times (\P^n)^{r+1}\mid
v_i\in\Lambda_{f(x)}\ ,\ i=0,\ldots,r\ \}$$ et notons $p:Z\ra X$ et
$q:Z\ra(\P^n)^{r+1}$ les deux projections. Comme $p$ est un fibr\'e en
produits de grassmanniennes, $Z$ est irr\'eductible; dans le cadre
topologique, on a aussi $\pi_1(Z)\dra\pi_1(X)$. On a
$$f^{-1}(\Sigma_{L_0}\cap\cdots\cap \Sigma_{L_r})=p\bigl(
q^{-1}(L_0\times\cdots\times L_r)\bigr)\ .$$ \ind Soit $I$ une partie de
$\{ 0,\ldots,r\}$ de cardinal $s+1\ge 1$;
 minorons la dimension de $p_Iq(Z)$. L'hypoth\`ese faite entra\^ine
$$f^{-1}(\Sigma_{M_0}\cap\cdots\cap\Sigma_{M_s})\ne\vide$$
 pour tous sous-espaces lin\'eaires $M_0,\ldots,M_s$ de $\P^n$ tels que
$\codim(M_i)=\codim (L_i)$ pour $i>0$ et $\codim(M_0)=\codim (L_0)+1$. On a
donc aussi $q^{-1}(M_0\cap\cdots\cap M_s)\ne\vide$
 et la prop. 3.1 entra\^ine  $$\dim
p_{\{0,\ldots,s\}}q(Z)\ge\sum_{i=0}^s\codim (M_i)\ .$$   \ind Puisque
$q(Z)$ est invariant par permutation des facteurs, on obtient $$\dim
p_Iq(Z)=\dim p_{\{0,\ldots,s\}}q(Z)\ge\sum_{i=0}^s\codim
(M_i)=1+\sum_{i=0}^s\ell_i >\sum_{i\in I}\ell_i\ .$$ \ind Le th. 1.3 et le
cor. 2.3 permettent de conclure dans ce cas.

\ind {\it Supposons maintenant $s>0$}. Pour $i=0,\ldots,s-1$, l'espace $L_i$
est r\'eduit \`a un point; soit $L$ l'espace lin\'eaire engendr\'e par
$L_0,\ldots,L_{s-1}$. Posons
$$\Sigma=\Sigma_{L_0}\cap\cdots\cap\Sigma_{L_{s-1}}=\{\ u\in  G(d,\P^n) \mid
\Lambda_u\supset L\ \}\ .$$ \ind Supposons les $L_i$ g\'en\'eraux; comme
$[\overline{f(X)}]\cdot \sigma_{n-d}^{s+1}\ne 0$, le th. 8.1 (qui n'est
dans ce cas que la version duale du th. 6.2) entra\^ine que
$X'=f^{-1}(\Sigma)$ est irr\'eductible. Si $L_i\cap L\ne\vide$, on a
$\Sigma_L\subset\Sigma_{L_i}$, de sorte que l'on peut supposer $L_i$
disjoint de $L$ pour $i=s,\ldots,r$. Soit $\P^{n-e}$ un sous-espace
lin\'eaire de $\P^n$ contenant $L_s,\ldots,L_r$ et disjoint de $L$. Il
existe un morphisme $f':X'\ra G(d-e,\P^{n-e})$ tel que
$\Lambda_{f'(x')}=\Lambda_{f(x)}\cap\P^{n-e}$, auquel il suffit d'appliquer
le cas d\'ej\`a trait\'e pour montrer a). On en d\'eduit b) et c) comme
d'habitude.\cqfd
 \bigskip {\bf 9. Conclusion}

\ind Soient $X$ et $Y$ des vari\'et\'es irr\'eductibles compl\`etes, et
$f:X\ra  G(d,\P^n)$ et $g:Y\ra  G(d,\P^n)$ des morphismes. Sous quelles
hypoth\`eses sur $f(X)$ et $g(Y)$ peut-on assurer que $X\times_{G(d,{\bf
P}^n)}Y$ est connexe?
 Je n'ai pas r\'eussi \`a \'enoncer une conjecture qui contienne tous les
r\'esultats de cet article. Lorsque $f$ est surjective, il suffit que $g$
soit non constante; mais en g\'en\'eral, m\^eme lorsque $f(X)$ est
encombrante, la condition $\dim f(X)+\dim g(Y)>\dim G(d,\P^n)$ ne suffit
pas. Qu'en est-il lorsque $f(X)$ et $g(Y)$ sont toutes deux encombrantes?
Le premier cas non trivial est celui o\`u $f(X)$ est un diviseur et $g(Y)$
une surface encombrante de $G(1,\P^3)$ (\cad\ distincte de $\Sigma_2$ et de
$\Sigma_{1,1}$).

\bigskip\bigskip \centerline{\pc R\'EF\'ERENCES}

\hangindent=1cm
 [D1] Debarre, O., {\it Th\'eor\`emes de connexit\'e et vari\'et\'es
ab\'eliennes},  \`a para\^\i tre dans Am. J. of Math (1995).

\hangindent=1cm
 [D2] Debarre, O., {\it On Subvarieties of Abelian Varieties}, \`a
para\^itre.

\hangindent=1cm
 [De] Deligne, P., {\it Le groupe fondamental du compl\'ement d'une courbe
plane n'ayant que des points doubles ordinaires est ab\'elien}, S\'eminaire
Bourbaki, Exp. \no 543, 1979/80, Springer Lecture Notes 842, 1981, 1--10.

\hangindent=1cm [F] Faltings, G.,  {\it Formale Geometrie und homogene
R\"aume}, Invent. Math. {\bf 64} (1981),\break 123--165.

\hangindent=1cm [Fu] Fulton, W.,  Intersection Theory, Springer Verlag,
Berlin, 1984.

\hangindent=1cm [FH] Fulton, W., Hansen, J., {\it A connectedness theorem
for projective varieties, with applications to intersections and
singularities of mappings}, Ann. of Math. {\bf 110} (1979), 159--166.

\hangindent=1cm [FL1] Fulton, W., Lazarsfeld, R., {\it Connectivity and its
Applications in Algebraic Geometry}, in  Algebraic Geometry,  Proceedings
of the Midwest Algebraic Geometry Conference, Chicago 1980, Springer
Lecture Notes 862, 1981, 26--92.

\hangindent=1cm [FL2] Fulton, W., Lazarsfeld, R., {\it Positive polynomials
for ample vector bundles}, Ann. of Math. {\bf 118} (1983), 35--60.

\hangindent=1cm [FMSS] Fulton, W., Mac Pherson, R., Sottile, F., Sturmfels,
B., {\it Intersection theory on spherical varieties},  J. Alg. Geom. {\bf
4} (1995), 181--193.

\hangindent=1cm [GL] Gaffney, T., Lazarsfeld, R., {\it On the Ramification
of Branched Coverings of $\P ^n$},  Invent. Math. {\bf 59} (1980), 53--58.

\hangindent=1cm [GH] Griffiths, P., Harris, J., Principles of Algebraic
Geometry, Wiley, New-York, 1978.

\hangindent=1cm [G] Gromov, M., {\it Convex Sets and K\"ahler Manifolds},
in Advances in Differential Geometry and Topology, World Scientific
Publishing, Teaneck, N.J., 1990.

\hangindent=1cm [Gr1] Grothendieck, {\it El\'ements de G\'eom\'etrie
Alg\'ebrique IV, 2},  Publ. Math. I.H.E.S. 24, 1965.

\hangindent=1cm [Gr2] Grothendieck, {\it El\'ements de G\'eom\'etrie
Alg\'ebrique IV, 3},  Publ. Math. I.H.E.S. 28, 1966.

\hangindent=1cm
 [H] Hansen, J., {\it A connectedness theorem for flagmanifolds and
Grassmannians}, Am. J. Math. {\bf 105} (1983), 633--639.

\hangindent=1cm [Ha] Harris, J.,  Algebraic Geometry,  Springer Verlag,
New-York, 1992.

\hangindent=1cm
 [HS] Hernandez, R., Sols, I., {\it Connectedness of intersections of
special Schubert varieties}, Manusc. Math. {\bf 83} (1994), 215--222.

\hangindent=1cm
 [Ho] Hovanski, A., {\it Fewnomials and Pfaff manifolds}, I.C.M. 1983,
Warszawa (1984), 549--565.

\hangindent=1cm [J]	Jouanolou, J.-P., {\it Th\'eor\`emes de Bertini et
applications},  Prog. Math. {\bf 42}, Birkh\"auser, 1983.

\hangindent=1cm [K]  Kleiman, S., {\it The transversality of a general
translate}, Comp. Math. {\bf 28} (1978), 287--297.

\hangindent=1cm [L] Lazarsfeld, R., Ph.D. thesis, Brown University, June
1980.

\hangindent=1cm [N]	Nori, M., {\it Zariski's conjecture and related
problems},  Ann. Sci. Ecole Norm. Sup. {\bf 16} (1983), 305--344.

\hangindent=1cm [PS] Paranjape, K.H., Srinivas, V., {\it Self maps of
homogeneous spaces}, Inv. Math. {\bf 98} (1989), 425--444.

\hangindent=1cm  [T]	Teissier, B., {\it Bonnesen-type inequalities in
algebraic geometry}, in  Seminar on Differential Geometry,  Princeton
University Press {\bf 102} (1982), 85--105.

 \bye